\documentclass[mathematics,article,submit,moreauthors,pdftex,10pt,a4paper]{mdpi}
\firstpage{1}
\articlenumber{x}
\doinum{10.3390/------}
\pubvolume{xx}
\pubyear{2018}
\copyrightyear{2018}
\history{Received: date; Accepted: date; Published: date}

\usepackage{amsmath,amsthm,amsopn,amstext,amscd,amsfonts,amssymb}
\usepackage{bm}
\newcommand{\vt}{\Delta_{i}^{j+1,j}}
\newcommand{\vti}{\Delta_{i}^{n_i,1}}
\def\build#1#2#3{\mathrel{\mathop{#1}\limits^{#2}_{#3}}}
\newcommand{\sumij}{\sum_{i=1}^d\sum_{j=1}^{n_i-1}}

\Title{A note on estimation of multi-sigmoidal Gompertz functions with random noise}

\Author{Patricia Román-Román $^{1,\ddagger}$\orcidA{}, Juan José Serrano-Pérez$^{2,\ddagger}$\orcidB{}, and Francisco Torres-Ruiz$^{1,\ddagger}$\orcidC{}}

\AuthorNames{Patricia Román-Román, Juan José Serrano-Pérez and Francisco Torres-Ruiz}

\address{%
$^{1}$ \quad Departamento de Estad{\'\i}stica e Investigaci{\'o}n Operativa, Universidad de Granada, Espa{\~n}a; Instituto de Matemáticas de la Universidad de Granada (IEMath-GR); proman@ugr.es; fdeasis@ugr.es

$^{2}$ \quad Departamento de Estad{\'\i}stica e Investigaci{\'o}n Operativa, Universidad de Granada, Espa{\~n}a; jjserra@ugr.es}

\corres{Correspondence: proman@ugr.es; Tel.: +34-958240491}

\secondnote{These authors contributed equally to this work.}

\abstract{The behaviour of many dynamic real phenomena shows different phases, with each one following a sigmoidal type pattern. This requires studying sigmoidal curves with more than one inflection point. In this work, a diffusion process is introduced whose mean function is a curve of this type, concretely a transformation of the well-known Gompertz model after introducing in its expression a polynomial term. The maximum likelihood estimation of the parameters of the model is studied, and various criteria are provided for the selection of the degree of the polynomial when real situations are addressed. Finally, some simulated examples are presented.}


\begin{document}

\section{Introduction}

Growth curves with sigmoidal behaviour are widely used for data analysis across several fields of application. Given the variety of approaches employed in the profuse literature dealing with this subject, several mathematical models have been proposed for their study.

In general terms, a sigmoidal function is a function defined in the real line, bounded and differentiable with positive derivative. Its graph has a typical \textit{S} shape showing a slow growth at the beginning, followed by a fast (exponential) growth that slows down gradually until it reaches an equilibrium value (usually named \textit{carrying capacity} or level of saturation). Often, the sigmoid function refers to the particular case of the logistic function $f(t)=k/(1+be^{-r\,t})$, although there is a great variety of curves with these characteristics. One of them is the Gompertz function, which is used in the modeling of systems that are saturated for large values of $t$, and whose most general expression is $f(t)=a\exp(-b\exp(-ct))$.

New sigmoidal curves have been introduced over the years, and their application has extended to several new fields. Among them we may mention hyperbolic curves (Eby et al. \cite{Eby10}), sigmoidal hyperbolic functions of the first and second kind (Menon et al. \cite{Men94}), or the beta growth function (Yin et al. \cite{Yin03}). Regarding their application to new fields, researchers have looked into the diffusion of innovations (Giovanis and Skiadas \cite{Gio07}); the calculation of oil production peaks (Gallagher \cite{Gal11}); predicting changes in language (Yokohama and Sanada \cite{Yok09}); analyzing fatigue and fractures in materials and structures (Paolino and Cavatorta \cite{Pao12}); etc.

Traditionally, most of the sigmoidal growth models cited above have arisen from the solution of ordinary differential equations. In that sense, they are deterministic and do not include other information than that provided by the variable under study. In order to incorporate these influences, the idea of growth in a random environment emerged (see Ricciardi, \cite{Ric79} and references therein). Thus, the so-called dynamic growth models appeared, and among them diffusion processes. Some of these diffusion models emerge as solutions to a stochastic differential equation after modifying a deterministic one by introducing in it a term of white noise. Other diffusion processes are constructed in such a way that their mean function is a certain sigmoidal growth curve. Usually, those of the first type retain the name of the corresponding deterministic equation of origin. For instance, Schurz \cite{Sch07}, collect a wide variety of logistic diffusion processes. However, for most of them the stochastic differential equation does not have an explicit solution. For this reason, Román and Torres \cite{Rom12} constructed a logistic-type stochastic differential equation in the second sense (its mean is a logistic function). The same situation has presented itself in the context of other growth curves and related diffusion processes. Such is the case of the Gompertz process, which was introduced by \cite{Cap74}. In such process, the upper limit of the curve is independent of the initial value of the population under study, something not always verified in real situations. For this reason, Gutiérrez et al. \cite{Gut07} introduced a new Gompertz-type process in which the \textit{carrying capacity} of the system depends on the initial state. This line of action has also been applied to the Bertalanffy curve (see Quiming et al. \cite{Qim07} and Román et al. \cite{Rom10}), the Hubbert curve (Luz-Sant'Ana et al. \cite{Ist17}), and, more recently, to the hyperbolastic curve of type I (Barrera et al. \cite{Bar18}).

Dynamic models are used in the fields in which the deterministic case has proved to be useful in fitting sigmoidal behaviour patterns to observed data. Researchers have developed several methods of estimation for these dynamic models. As far as maximum likelihood estimation methods are concerned, we can differentiate between those that take as a starting point the stochastic differential equation related to the model (usually known as continuous sampling methods) from those who build the likelihood function from the transition density functions of the process (discrete sampling methods). Alternatively, some authors have dealt with inference from a Bayesian perspective (Tang and Heron \cite{Tan08}).

An interesting aspect in this type of processes is the possibility of introducing, into their infinitesimal moments, time functions that allow us to regulate the evolution of the variable under study. Given that the functional form of such functions is not known, several strategies have been devised for their estimation. Some work carried out along this line includes studies by Albano et al. \cite{Alb11,Alb15} and Román et al. \cite{Rom16} which centered on modifications of the Gompertz process.

There are multiple real situations in which the maximum level of growth is reached after successive stages, in each of which there is a deceleration followed by an explosion of the exponential type.
For this reason, the use of sigmoidal curves with more than one inflection point is a good approach. A typical example of this behaviour is observed in the growth of various fruit species, such as stone fruits (Álvarez and Boché \cite{Alv99}). Cairns et al. \cite{Cai08} used double-sigmoidal models to study fatigue profiles in mouse muscles, while Amorim et al. \cite{Amo93} detected this type of behaviour in the different phases in which the fungus Ustilago Scitaminea Sydow infects the sugarcane and produces its characteristic smut.

The way that this multi-sigmoidal behaviour is modeled is far from unique. For example, Roper \cite{Rop00} used hyperbolic functions to study the transition between various temperature states in certain geological zones. Other authors have achieved this goal by including terms that define inflection points and additional parameters (Lipovetsky \cite{Lip10}). However, these models have not addressed the incorporation of external information to the variable under study, similarly to the dynamic models already mentioned. In this paper we address this problem by introducing a diffusion process whose mean obeys a pattern of multi-sigmoidal behaviour. In particular, we will deal with the case of Gompertz growth with multiple inflection points, following the idea mentioned in \cite{Amo93} for the case of the generalized monomolecular and Gompertz curves.

The rest of the paper is organized as follows: in Section 2 the multi-sigmoidal Gompertz curve is introduced by including a polynomial in the usual expression of the curve. In Section 3, the Gompertz multi-sigmoidal diffusion process is defined. To this end, the lognormal diffusion process with exogenous factors is considered, since it allows us to model behavioral patterns that verify the properties exhibited by the curve. The estimation of the process, which is performed by maximum likelihood using discrete sampling, is the subject of Section 4. The matter of obtaining initial solutions to solve the resulting system of equations deserves special attention. Other important aspect is to determine the degree of polynomial that should be considered since, in general, this aspect will be unknown in real applications. To this end, some criteria are considered in that section. Finally, in Section 5, some simulation examples are considered.

\section{Multi-sigmoidal Gompertz curve}

Let $Q_{\bm\beta}(t)=\sum_{l=1}^{p}\beta_lt^l$ be a $p$-degree polynomial , where
$\bm \beta=(\beta_1,\ldots,\beta_p)^T$ ($p>1$) denotes a real parametric vector with positive leading coefficient $\beta_p$. We define the multi-sigmoidal Gompertz function as
\begin{equation}
\label{curva}
f_{\bm\theta}(t)=k\exp\left(-\alpha\,e^{-Q_{\bm\beta}(t)}\right), \qquad  t\geq t_0\geq 0, \qquad \alpha,k>0, \qquad \bm{\theta}=(\alpha,\bm\beta^T)^T.
\end{equation}

Denoting $P_{\bm\beta}(t)=\frac{dQ_{\bm\beta}(t)}{dt}$, curve (\ref{curva}) satisfies the ordinary linear differential equation
\begin{equation}
\label{EDO1}
\frac{df_{\bm\theta}(t)}{dt}=f_{\bm\theta}(t)h_{\bm\theta}(t),
\end{equation}
where
\begin{equation}
\label{Funcionh}
h_{\bm\theta}(t)=\alpha\,P_{\bm\beta}(t)\,e^{-Q_{\bm\beta}(t)}.
\end{equation}

Taking into account that $\ln f_{\bm\theta}(t)=\ln k-\alpha\,e^{-Q_{\bm\beta}(t)}$, the above equation can be expressed as
\begin{equation}
\label{EDO2}
\frac{df_{\bm\theta}(t)}{dt}=f_{\bm\theta}(t)\left(\ln k-\ln f_{\bm\theta}(t)\right)P_{\bm\beta}(t).
\end{equation}

However, the resolution of both equations, with initial condition $f_{\bm\theta}(t_0)=f_0>0$, leads to two expressions of curve \eqref{curva}. Indeed, the solution of \eqref{EDO1} is
\begin{equation}
\label{Curva1}
f_{\bm\theta}(t)=f_0\exp\left(-\alpha\,\left(e^{-Q_{\bm\beta}(t)}-e^{-Q_{\bm\beta}(t_0)}\right)\right),
\end{equation}
while the solution of \eqref{EDO2} is
\begin{equation}
\label{Curva2}
f_{\bm\theta}(t)=\exp\left(\ln\,k\left(1-e^{-Q_{\bm\beta}(t)-Q_{\bm\beta}(t_0)}\right)+\ln f_0\,e^{-Q_{\bm\beta}(t)-Q_{\bm\beta}(t_0)}\right).
\end{equation}

Equation \eqref{EDO2} is a generalization of the classical gompertzian differential equation, giving rise to the curve \eqref{Curva2}, used by authors like Ricciardi et al. \cite{Ric83} in the case $Q_\beta(t)=\beta t$. On the other hand, equation \eqref{EDO1} is a linear differential equation of the Malthusian type whose solution generalizes the expression of the Gompertz curve used by authors such as Laird \cite{Lai65} and Gutiérrez et al. \cite{Gut07}.

The main difference between \eqref{Curva1} and \eqref {Curva2} lies in their limit value, which in the first case is $k(\bm\theta)=f_0\exp\left(\alpha\, e^{-Q_{\bm\beta}(t_0)}\right) $, and $k$ in the second. This may lead to the choice of either expression depending on the knowledge available about the influence of initial value $f_0$ on the limit value.  This is the case when the phenomenon under study shows Gompertz-type growth and several sample paths are available, each with a common growth pattern but with different initial values and a different limit value (for example, the particular weight of each individual of the same species). In the rest of the paper, we will consider the situation in which the \text{carrying capacity} of the system modeled by the curve depends on the initial value of the population under study. So, \eqref{Curva1} can be expressed as $f_{\bm\theta}(t)=k(\bm\theta)g_{\bm\theta}(t)$ where $g_{\bm\theta}(t)=
\exp\left(-\alpha\, e^{-Q_{\bm\beta}(t)}\right)$, verifying $\lim_{t\rightarrow\infty}g_{\bm\theta}(t)= 1$.

Focusing on the general expression given by \eqref{curva}, from \eqref{EDO1} it follows that the growth intervals of the curve depend on the roots of the equation $P_{\bm\beta}(t)=0$. As for the inflection points, the candidates will be the solutions of $\frac{d^2f_{\bm\theta}(t)}{dt^2}=0$, which results in solving the equation
\begin{equation}
\label{DerivadaSegunda}
  \frac{dP_{\bm\beta}(t)}{dt}=P^2_{\bm\beta}(t)\left(1-\alpha\,e^{-Q_{\bm\beta}(t)}\right).
\end{equation}

\begin{figure}[h!]
$$
\begin{array}{ccc}
& a) & \\
\includegraphics[height=5cm, width=6cm]{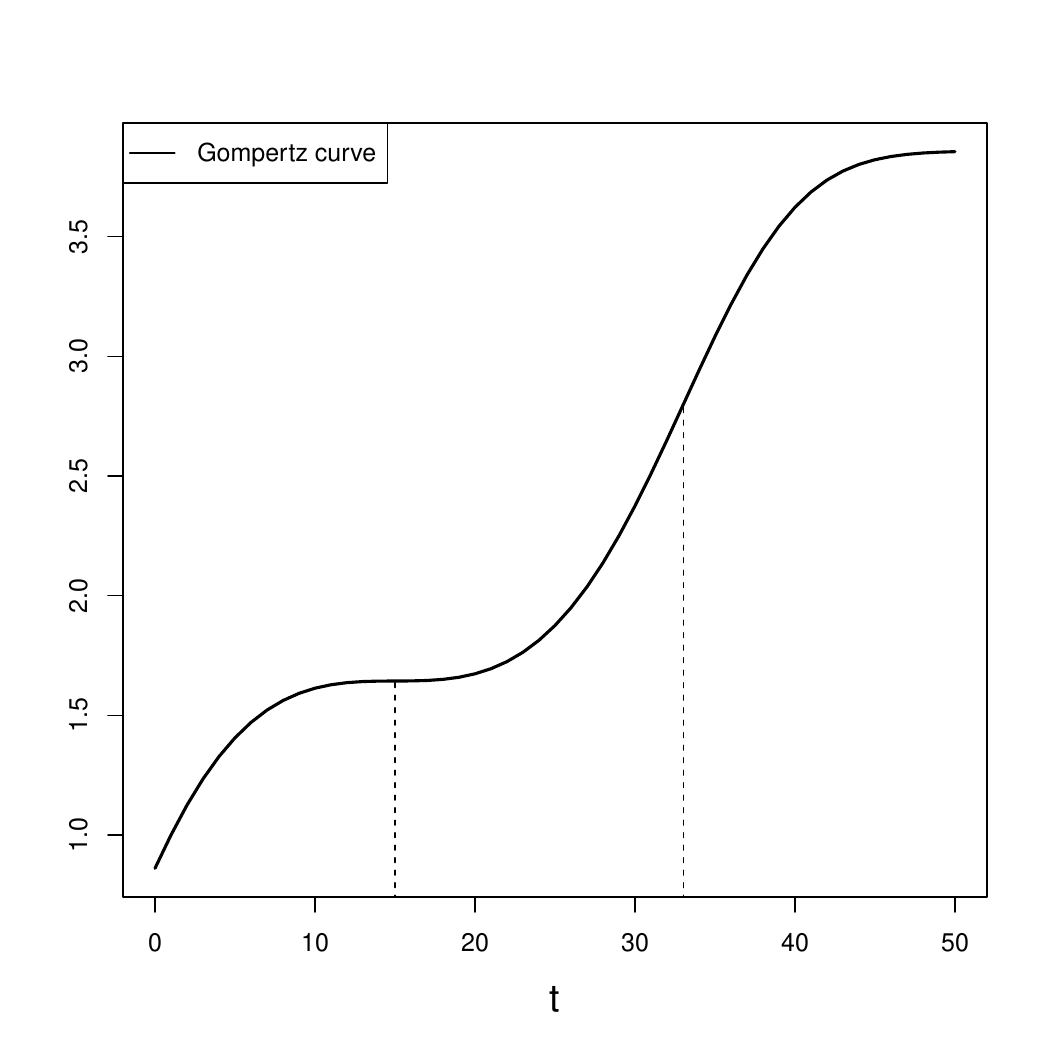} & & \includegraphics[height=5cm, width=6cm]{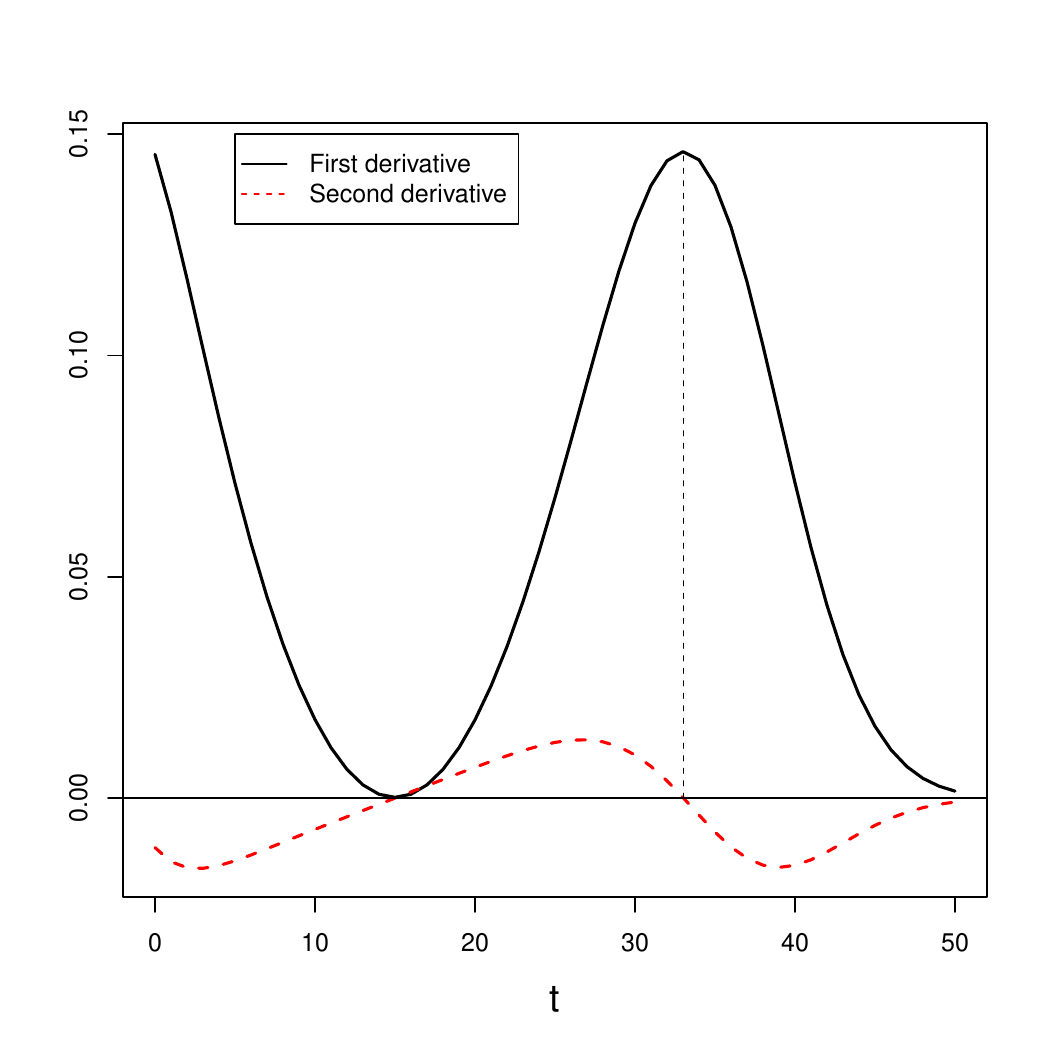}\\
& b) & \\
\includegraphics[height=5cm, width=6cm]{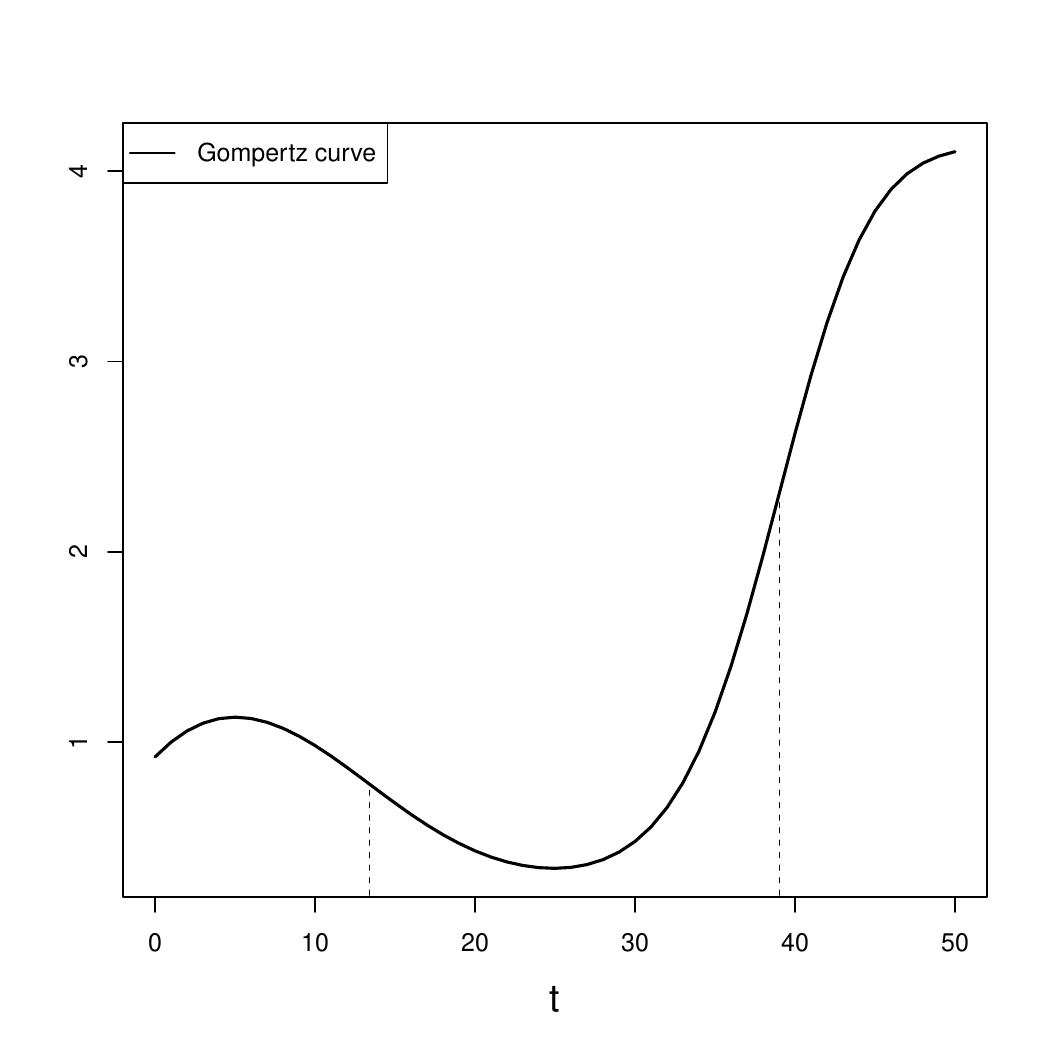} & & \includegraphics[height=5cm, width=6cm]{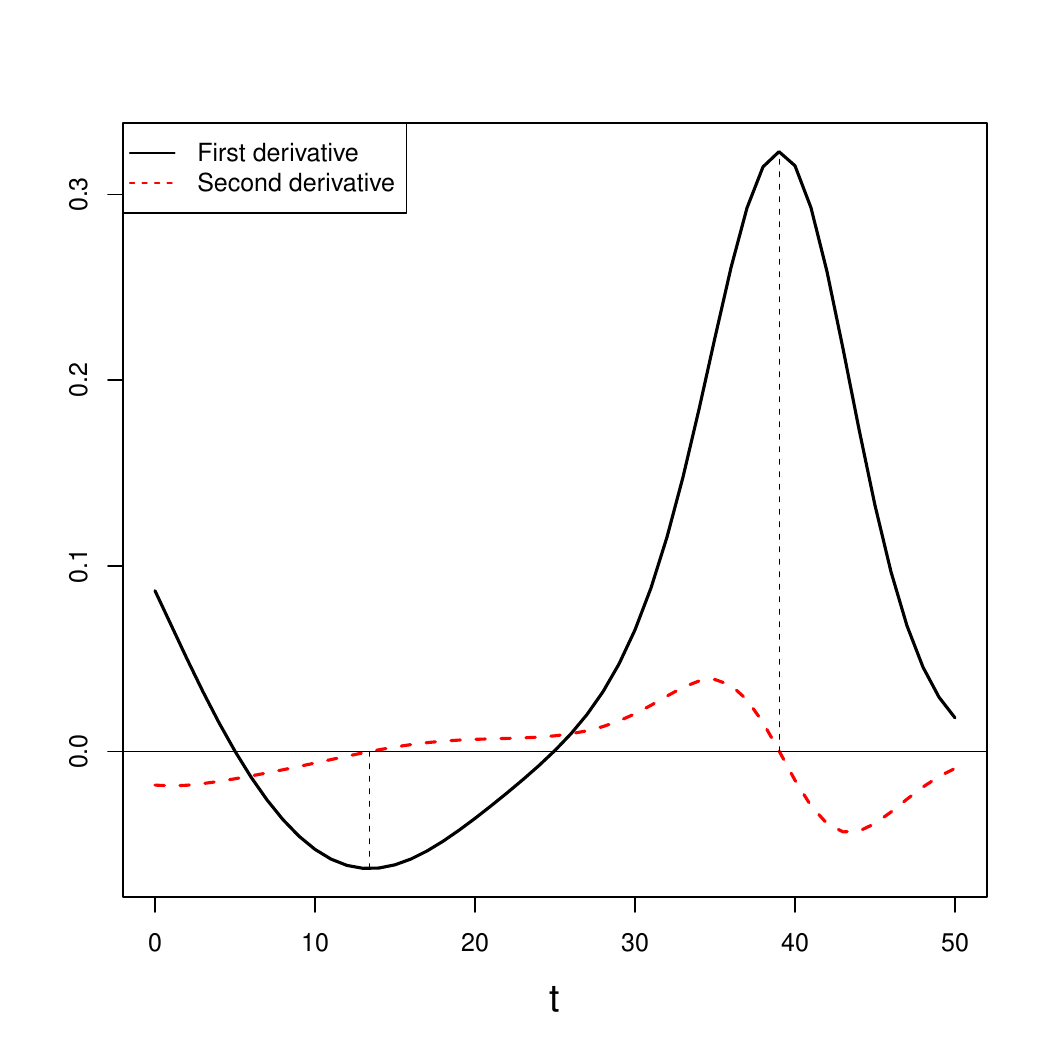}\\
& c) & \\
\includegraphics[height=5cm, width=6cm]{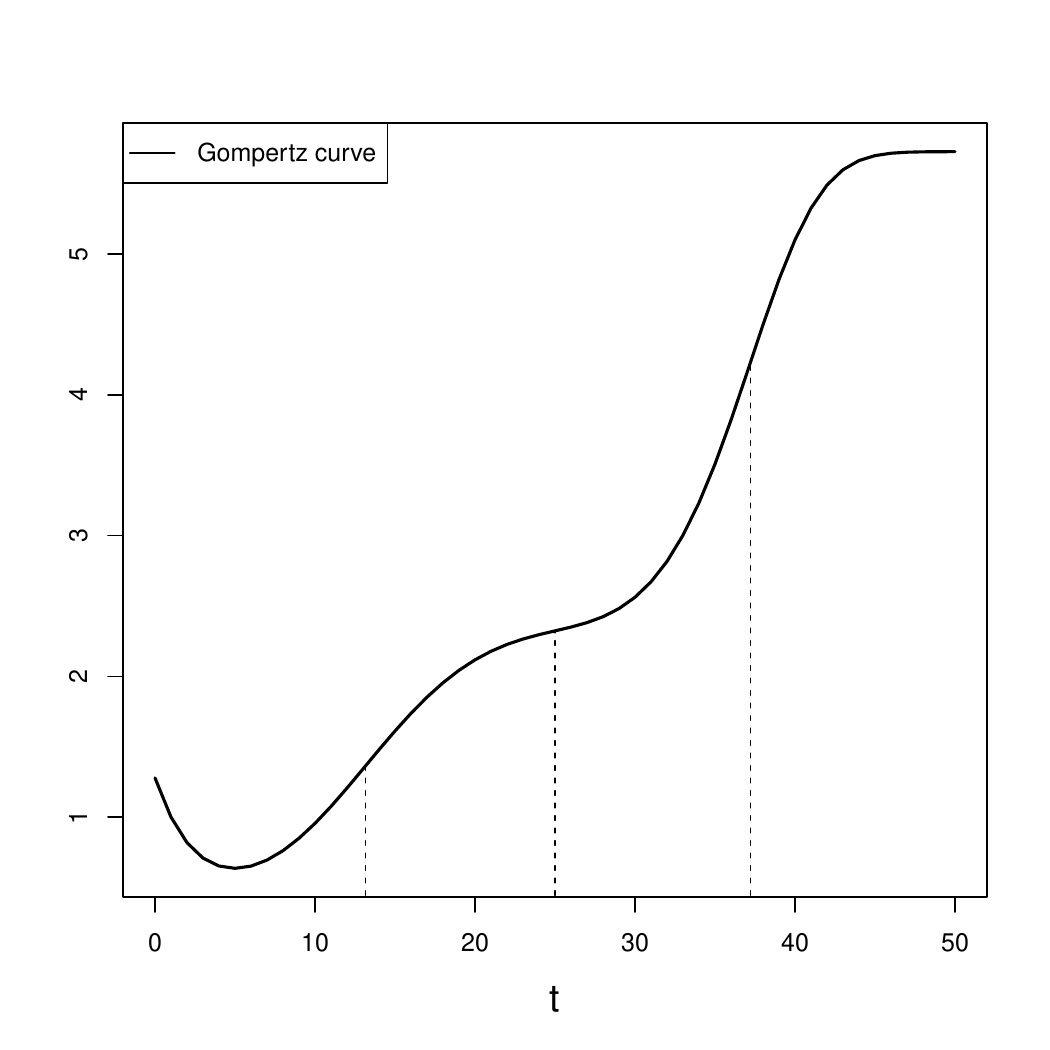} & & \includegraphics[height=5cm, width=6cm]{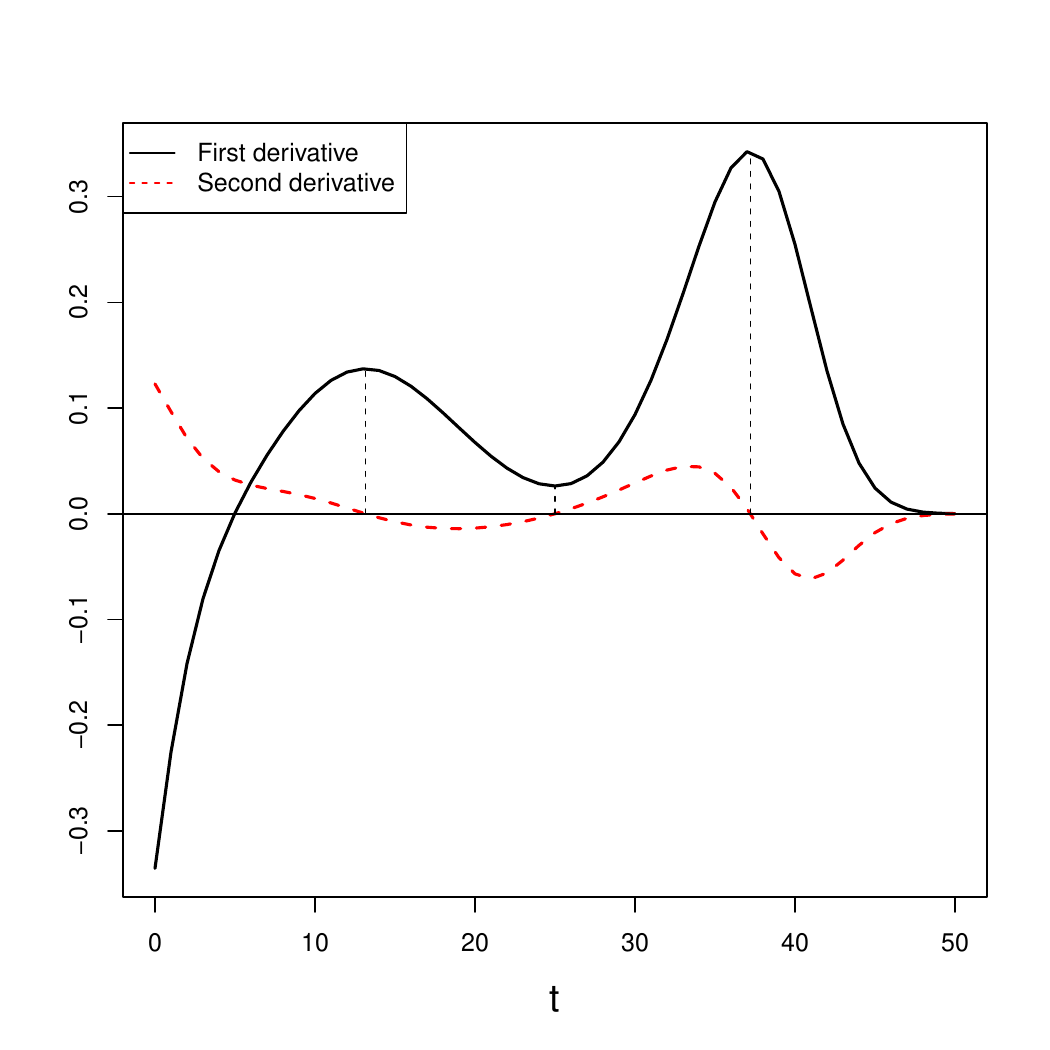}\\
\end{array}
$$
\caption{\small Several examples of multi-sigmoidal Gompertz curves.}
\label{examples}
\end{figure}

Figure \ref{examples} shows some possible situations for various choices of the polynomial $Q_{\bm\beta}(t)$. In each case, the Gompertz curve is represented together with its first and second derivatives. In particular, figure a) represents the case in which $P_{\bm\beta}(t) $ has no roots and the curve is strictly increasing and presenting two inflection points, as in case b), although in the latter case the curve presents both decreasing and increasing intervals. Finally, figure c) shows an example with three inflection points.

\section{Multi-sigmoidal Gompertz diffusion process}

Assuming that \eqref{DerivadaSegunda} has at least one solution where $f_{\bm\theta}$ presents an inflection point, this function is found in the conditions listed in Román and Torres \cite{Rom15}, which ensures that the growth phenomenon represented by the curve can be modeled by a non-homogeneous lognormal diffusion process whose mean function is such a function. In Román et al. \cite{Rom18}, a general study of this process is carried out, including the distribution and main characteristics as well as aspects related to inference. Following the notation used in that paper, we define the multi-sigmoidal Gompertz process as a diffusion process $\{X(t); t\in I\} $ that takes values in $(0,+\infty)$ and with infinitesimal moments
\begin{equation*}
\begin{array}{l}
A_1(x,t)=h_{\bm\theta}(t)x\\
A_2(x)={\sigma}^{2}x^{2}, \qquad \sigma>0
\end{array}
\end{equation*}
where $I=[t_0, +\infty)$ is a real interval ($ t_0 \geq 0 $) and $\Theta\subseteq\mathbb{R}^{p+1}$ is an open set such that $\bm\theta=(\alpha, \bm\beta^T)^T \in \Theta$, where $h_ {\bm\theta} $ is given by \eqref{Funcionh}. This process is determined from the stochastic differential equation
$$
dX(t)=h_{\bm\theta}(t)X(t) dt+ \sigma X(t) dW(t), \qquad X(t_0)=X_0,
$$
where $W(t)$ is a Wiener process (Brownian motion), independent of the initial condition $X_0=X(t_0)$, $t\geq t_0$. The solution to this equation can be expressed as
\begin{equation}
\label{SolEDE}
X(t)=X_0 \, \exp \left(H_{\bm\xi}(t_0,t)+\sigma(W(t)-W(t_0)) \right), \qquad t\geq t_0
\end{equation}
where
$$
H_{\bm\xi}(s,t)=\int_{s}^{t}h_{\bm\theta}(u) d u-\frac{\sigma^2}{2}(t-s)=-\alpha\left(e^{-Q_{\bm\beta}(t)}-e^{-Q_{\bm\beta}(s)}\right)-\frac{\sigma^2}{2}(t-s), \qquad s<t, \,\bm\xi=(\bm\theta^T,\sigma^2)^T.
$$

Regarding the distribution of the process, if $X_0 $ is distributed according to a lognormal distribution $\Lambda_1 \left[\mu_0; \sigma_0^2\right]$, or $X_0$ is a degenerate variable (i.e. $ P[X_0=x_0]=1$),
the finite dimensional distributions of the process are lognormal. Thus, $\forall n\in\mathbb{N}$ and $t_1<\cdots<t_n$, vector $(X(t_1),\ldots, X(t_n))^T $ is distributed according to an $n$-dimensional lognormal distribution $\Lambda_n[\bm\varepsilon, \bm\Sigma]$, where the components of the vector $\bm\varepsilon$ are $\varepsilon_i=\mu_0+H_{\bm\xi}(t_0,t_i) $, $i=1,\ldots, n$, being
$\sigma_{ij}=\sigma_0^2+\sigma^2(\min(t_i, t_j)-t_0)$, $i,j=1, \ldots, n$, those of the matrix $\bm\Sigma$.

From the two-dimensional distributions $(X(s),X(t))^T$, $s<t$, the transitions of the process can be obtained, which are also lognormal; concretely,
\begin{equation}
\label{Transicion}
X(t) \mid X(s) = y \rightsquigarrow \Lambda_1 \left( \ln y + H_{\bm\xi}(s,t), \sigma^2 (t-s)\right), \qquad s<t.
\end{equation}

Once the distribution of the process has been established, different characteristics associated with it can be calculated, including the mean and conditioned mean functions, whose expressions are
\begin{equation*}\label{media}
  m(t)=E[X(t)]=E[X_0]\exp\left(-\alpha\,\left(e^{-Q_{\bm\beta}(t)}-e^{-Q_{\bm\beta}(t_0)}\right)\right)
\end{equation*}
and
\begin{equation*}\label{media}
  m(t|t_0)=E[X(t)|X(t_0)=x_0]=x_0\exp\left(-\alpha\,\left(e^{-Q_{\bm\beta}(t)}-e^{-Q_{\bm\beta}(t_0)}\right)\right)
\end{equation*}
respectively, being able to verify that both functions are of the type introduced in the previous section. The expression that adopts other characteristics can be consulted in Román et al. \cite{Rom18}.

\section{Maximum likelihood estimation}

In this section we will deal with the estimation of the parameters of the process by maximum likelihood. In Román and Torres \cite{Rom18}, a general treatment of this question is carried out for the non-homogeneous lognormal process, of which the process dealt with in this paper is a particular case. Below we summarize the main results obtained in the aforementioned paper, adapting them to the actual case.

The starting point is a discrete sampling of the process based on $d$ sample paths observed at time instants $t_{ij}$, $(i=1, \ldots, d, \, j=1, \ldots, n_i)$. Note that the observation time instants do not have to be the same in each trajectory, although we will suppose that $t_{i1}=t_0$, $i=1, \ldots, d$. Denote by $\mathbf{X}=\left(\mathbf{X}_1^T | \cdots | \mathbf{X}_d^T \right)^T$, where $\mathbf{X}_i^T$ is the vector that contains the variables of the $i$-th sample-path, that is
$\mathbf{X}_i=(X(t_{i1}),\ldots,X(t_{i,n_i}))^T$, $i=1,\ldots,d$.

Assuming that the distribution of $X(t_0)$ is lognormal $\Lambda_1 (\mu_1, \sigma_1^2)$, and taking into account the transitions of the process, the probability density function of $\mathbf{X}$ is given by
\begin{align*}
f_\mathbf{X}(\mathbf{x})&=\prod_{i=1}^{d}\displaystyle\frac{\exp\left(-\frac{[\ln x_{i1}-\mu_1]^2}{2\sigma_1^2}\right)}{x_{i1}\sigma_1\sqrt{2\pi}}
\prod_{j=1}^{n_i-1}\displaystyle\frac{\exp\left(-\frac{\left[\ln\left(x_{i,j+1}/x_{ij}\right)-m_{\bm\xi}^{i,j+1,j}\right]^2}{2\sigma^2\vt}\right)}{x_{ij}\sigma\sqrt{2\pi\vt}},
\end{align*}
where $m_{\bm\xi}^{i,j+1,j}$ and $\vt$ are given by
\begin{align*}\label{Nota1}
  m_{\bm\xi}^{i,m,n} & = H_{\bm\xi}(t_{in},t_{im})=-\alpha\,\phi_{i,m,n}^{\bm\beta,0}-\displaystyle\frac{\sigma^2}{2}\Delta_i^{m,n}\\
  \Delta_i^{m,n} &= t_{im}-t_{in}
\end{align*}
with
\begin{equation*}\label{phi}
  \phi_{i,m,n}^{\bm\beta,l}=t_{im}^l\,e^{-Q_{\bm\beta}(t_{im})}-t_{in}^l\,e^{-Q_{\bm\beta}(t_{in})}, \qquad l=0,1,...; \qquad i=1,\ldots,d; \qquad m,n\in\{1,\ldots,n_{i-1}\}, m>n.
\end{equation*}

Next, we consider the change of variables
\begin{align*}
\label{CambioVariable}
V_{0i}&=X_{i1},\,\, i=1,\ldots,d\\ \nonumber
V_{ij}&=(\vt)^{-1/2}\ln\frac{X_{i,j+1}}{X_{ij}},\,\, i=1,\ldots,d; j=1,\ldots,n_i-1
\end{align*}
that transforms vector $\mathbf{X}$ in $\mathbf{V}=\left[\mathbf{V}_0^T|\mathbf{V}_1^T|\cdots|\mathbf{V}_d^T\right]^T=\left[\mathbf{V}_0^T|\mathbf{V}_{(1)}^T\right]^T$, whit probability density function
\begin{equation}\label{DensidadV}
f_\mathbf{V}(\mathbf{v})=
\displaystyle\frac{\exp\left(-\frac{1}{2\sigma_1^2}
(\ln\mathbf{v}_0-\mu_1\mathbf{1}_d)^T(\ln\mathbf{v}_0-\mu_1\mathbf{1}_d)\right)}
{\displaystyle\prod_{i=1}^{d}v_{0i}\left(2\pi\sigma_1^2\right)^{\frac{d}{2}}}
\displaystyle
\frac{\exp\left(-\frac{1}{2\sigma^2}\left(\mathbf{v}_{(1)}-\bm\gamma^{\bm\xi}\right)^T\left(\mathbf{v}_{(1)}-\bm\gamma^{\bm\xi}\right)\right)}
{\left(2\pi\sigma^2\right)^{\frac{n}{2}}}
\end{equation}
where $\ln\mathbf{v}_0=(\ln v_{01},\ldots,\ln v_{0d})^T$, $n=\sum_{i=1}^d(n_i-1)$, $\mathbf{1}_d=(1,\ldots,1)^T_{d\times 1}$, and $\bm\gamma^{\bm\xi}$ is an $n$-dimensional vector with components $\gamma^{\bm\xi}_{ij}=(\vt)^{-1/2}m_{\xi}^{i,j,j+1}$, $i=1,\ldots,d; j=1,\ldots,n_i-1$.

For a fixed value $\mathbf{v}$, expression \eqref{DensidadV} provides the likelihood function, whose logarithm is
\begin{align*}
\label{Vero2}
 L_{\mathbf{v}}(\bm\eta,\bm\xi)&=-\frac{(n+d)\ln(2\pi)}{2}-\frac{d\ln\sigma_1^2}{2}-\sum_{i=1}^d\ln v_{0i}-\frac{\displaystyle\sum_{i=1}^d\left[\ln v_{0i}-\mu_1\right]^2}{2\sigma_1^2}
     -\frac{n\ln\sigma^2}{2}-\frac{Z_1+\Phi_{\bm\xi}-2\Gamma_{\bm\xi}}{2\sigma^2}
\end{align*}
where $\bm\eta=(\mu_1, \sigma_1^2)^T$ is the vector that contains the parameters of the initial distribution, being
$$
Z_1=\displaystyle\sumij v_{ij}^2, \qquad \Phi_{\bm\xi}=\displaystyle\sumij \frac{\left(m_{\bm\xi}^{i,j+1,j}\right)^2}{\vt},\qquad
\Gamma_{\bm\xi}=\displaystyle\sumij \frac{v_{ij}m_{\bm\xi}^{i,j+1,j}}{(\vt)^{1/2}}.
$$

Assuming that $\bm\eta$ and $\bm\xi$ are functionally independent, the estimate of $\bm\eta$ leads to
\[
  \widehat{\mu}_1=\displaystyle\frac{1}{d}\displaystyle\sum_{i=1}^{d}\ln v_{0i} \ \ \mbox{and} \ \
  \widehat{\sigma}_1^2=\displaystyle\frac{1}{d}\displaystyle\sum_{i=1}^{d}(\ln v_{0i}-\widehat{\mu}_1)^2,
\]
while that of $\bm\xi$ is obtained (see \cite{Rom18} for details) from the system of equations
\begin{align}
\label{eq1sistema}
&\Psi_{\bm\theta}-\Omega_{\bm\xi}=0\\
\label{eq2sistema}
&Z_1+\Phi_{\bm\xi}-2\Gamma_{\bm\xi}-\sigma^2Z_2+\sigma^2\Upsilon_{\bm\xi}=n\sigma^2
\end{align}
where
\begin{equation*}
\label{funciones}
  \Omega_{\bm\xi}=\frac{1}{2}\frac{\partial\Phi_{\bm\xi}}{\partial\bm\theta^T},\qquad
  \Psi_{\bm\theta}=\frac{1}{2}\frac{\partial\Gamma_{\bm\xi}}{\partial\bm\theta^T},\qquad
  \Upsilon_{\bm\xi}=-\frac{\partial\Phi_{\bm\xi}}{\partial\sigma^2},\qquad Z_2=-2\frac{\partial\Gamma_{\bm\xi}}{\partial\sigma^2}
\end{equation*}

Taking into account that
$$
\frac{\partial m_{\bm\xi}^{i,j+1,j}}{\partial\bm\theta^T}=\left(-\phi_{i,j+1,j}^{\bm\beta,0},\alpha\phi_{i,j+1,j}^{\bm\beta,1},\alpha\phi_{i,j+1,j}^{\bm\beta,p}\right),
$$
and the previous expressions of $\Omega_{\bm\xi}$ and $\Psi_{\bm\theta}$, the subsystem of equations \eqref{eq1sistema} remains in the form
\begin{equation}
\label{eq1sistema2}
  X_l^{\bm\beta}+\alpha Y_l^{\bm\beta}+\frac{\sigma^2}{2}W_l^{\bm\beta}=0, \qquad l=0,1,\ldots,p
\end{equation}
where, for $l=0,1,\ldots, p$ one has
$$
X_l^{\bm\beta}=\displaystyle\sumij \frac{v_{ij}\phi_{i,j+1,j}^{\bm\beta,l}}{(\vt)^{1/2} }, \qquad Y_l^{\bm\beta}=\displaystyle\sumij \frac{\phi_{i,j+1,j}^{\bm\beta,0}\phi_{i,j+1,j}^{\bm\beta,l}}{\vt},
\qquad W_l^{\bm\beta}=\displaystyle\sum_{i=1}^d\phi_{i,n_i,1}^{\bm\beta,l}.
$$

On the other hand, and since
$$
\Phi_{\bm\xi}  = \alpha^2Y_0^{\bm\beta}+\frac{\sigma^4}{4}Z_3+\alpha\sigma^2W_0^{\bm\beta}, \qquad
\Gamma_{\bm\xi}  =-\alpha\,X_0^{\bm\beta}-\frac{\sigma^2}{2}Z_2, \qquad
\Upsilon_{\bm\xi}  =-\alpha\,W_0^{\bm\beta}-\frac{\sigma^2}{2}Z_3,\qquad Z_3=\displaystyle\sum_{i=1}^d\vti
$$
equation \eqref{eq2sistema} transforms into

\begin{equation}
\label{Ecuacionalfa}
\sigma^2\left[n+\sigma^2Z_3/4\right]-\alpha\left[2X_0^{\bm\beta}+\alpha\,Y_0^{\bm\beta}\right]-Z_1=0.
\end{equation}

\subsection{System of equations and numerical computations}

The system of equations \eqref{eq1sistema2}-\eqref{Ecuacionalfa} can not be solved explicitly, and it is therefore necessary to use numerical methods, such as Newton-Raphson, for which an initial solution is required. Next, we present a strategy to achieve this, based on the information provided by the sample data.

As a matter of fact, and taking into account that the mean function of the process is a Gompertz multi-sigmoidal curve, as well as the expression \eqref{curva}, it follows
$$
\ln\ln \frac{k}{f_{\bm\theta}(t)}=\ln\alpha-Q_{\bm\beta}(t).
$$
Noting $m_i$ the values of the mean of the sample paths at $t_i$, we propose to fit, by linear regression, a polynomial taking as data the pairs of values $(t_i, \ln(\ln(k/m_i)))$. The estimated coefficients will provide the initial values for $\alpha$ and $\bm\beta$. Regarding $\sigma$, its initial estimation is based on the fact that for a lognormal distribution $\Lambda_1[\eta,\delta]$, the quotient between the arithmetic mean and the geometric one provides an estimation of $\delta$; concretely $\widehat{\delta}=2\ln(E[X]/E_G[X])$. Applying this result to the distribution of $X(t)$ we obtain, for each $t_i$, an estimate of $\sigma^2\,(t_i-t_0)$; that is, $\sigma^2_i=2\log (m_i/m^g_i)$, $i=1,\ldots$, where $m^g_i$ are the values of the geometric sample mean. Finally, the initial value of $\sigma$ is calculated by performing a simple linear regression of the $\sigma_i$ values against $t_i$.

In this procedure there are several questions that must be taken into account:
\begin{itemize}
  \item The value of $k$, in general, will not be known. Therefore, we suggest taking as an approximation the last value of the mean. However, it is possible that in real cases, and due to the fluctuations of the process, there could be values $m_i$ verifying $k\leq m_i$, so transformation $\ln(\ln(k/m_i))$ would not be determined. In such cases, usually not many in practical cases, those points must be removed from the regression analysis.
  \item Since, generally, the degree of the polynomial will not be known a priori, it is necessary to have some mechanism that will allow for its selection. To this end we propose a forward procedure, introducing polynomials in a consecutive way in the model. Each time a polynomial is introduced, a measure of the adjustment made is calculated and compared with the previous ones. If the adjustment is improved, the procedure continues; otherwise it stops. However, and even in this case, it is convenient to perform one more iteration due to the parity of the polynomial.
  \item Regarding the measures that can be used to evaluate the adjustment, we propose the following:
      \begin{itemize}
        \item The absolute relative errors between the sample mean of the process and the fitted mean for each estimated model
            $$
            RAE_j= \displaystyle\frac{1}{N}\sum_{i=1}^N\displaystyle\frac{|m_i-\widehat{E}[X^{(j)}(t_i)]|}{m_i},
            j=1,2,\ldots
            $$
        \item The resistor-average distance (see Johnson and Sinanovic \cite{Joh01}). This is a distance based on the Kullback–Leibler divergence, which will be used for calculating the distance between the sample distribution (available from the data) and that obtained from each estimated model. The expression for this measure is
            \begin{equation*}
            \label{DRA}
            D_{RA}(f_s||f_i)=\frac{D_{KL}(f_s||f_i)\,\cdot\,D_{KL}(f_i||f_s)}{D_{KL}(f_s||f_i)+D_{KL}(f_i||f_s)},
            \end{equation*}
            where $D_{KL}(f_s||f_i)$ denotes the Kullback–Leibler divergence between the sample distribution ($f_s$) and that for the $i$-th estimated model ($f_i$). Its expression is given by
            \begin{equation*}
            D_{KL}(f_s||f_i)=\frac{1}{2}\left[\log\left(\frac{\widehat{\sigma}^2(t_i-t_0)}{\sigma_i^2}\right)+
            \frac{\sigma_i^2}{\widehat{\sigma}^2(t_i-t_0)}+
            \frac{\left(\log m^g_i-\log \widehat{E[X_0]}-H_{\widehat{\bm\xi}}(t_0,t_i)\right)^2}{\widehat{\sigma}^2(t_i-t_0)}-1\right].
            \end{equation*}

            In practice, this is the expression of the distance that should be used since the theoretical model will not be known in real applications. However, for simulation studies the distance between theoretical and estimated models could be considered. In this case, the previous expression would be slightly modified.
        \item The Akaike information criterion (AIC) and the Bayes information criterion (BIC).
      \end{itemize}
\end{itemize}

\subsection{About $t_0$}

  Another interesting aspect to consider is that of the time instants, especially when they take high values. This can mainly affect the obtention of the initial values for the parameters since a polynomial regression has been proposed. One option is to apply orthogonal polynomials, as it is usual when considering this type of regression.

  One alternative is to consider a new diffusion process $\{Y(t); t\geq 0\}$ obtained from $\{X(t); t\geq t_0\}$ by considering a shift of length $t_0$ in time, that is $Y(t)=X(t+t_0)$, so that the original data can be considered as observations of the new process with an initial instant equal to zero. Let's see how this change affects the infinitesimal moments of the processes:

  In general, let $\{X(t); t\geq t_0\}$ be the original process and $\{Y(t); t\geq 0\}$ verifying $X(t)=Y(t-t_0)$. Denote by $A_m^X(x,t)$ and $A^Y_m(x,t)$ their respective $m$-th order infinitesimal moments. Taking into account the definition of infinitesimal moment of order $m$,
  $$
  A_m^X(x,t)=\build{\lim}{}{h\rightarrow 0} \displaystyle\frac{1}{h}  E\left[(X(t+h)-X(t))^m|X(t)=x\right],
  $$
  it verifies that $A_m^X(x,t)=A_m^Y(x,t-t_0)$. Obviously, the strategy of considering this translation over time will be useful when the resulting process is of the same type as the original.

  In the case of the multi-sigmoidal Gompertz process, let us consider $\{Y(t); t\geq 0\}$ with infinitesimal moments
  $$
  A_1^Y(x,t)=\eta\widetilde{P}_{\bm\gamma}(t)\,e^{-\widetilde{Q}_{\bm\gamma}(t)}, \qquad A_2^Y(x)=\sigma^2\,x^2
  $$
  where $\bm\gamma=(\gamma_1,\ldots,\gamma_p)^T$.

  Taking into account that
  $$
  \widetilde{Q}_{\bm\gamma}(t-t_0)=\sum_{l=1}^{p}\gamma_l\,(t-t_0)^l=\beta_0+\sum_{m=1}^{p}\beta_m\,t^m=
  \beta_0+Q_{\bm\beta}(t),
  $$
  where $\bm\gamma=(\gamma_1,\ldots,\gamma_m)$, being
  $$
  \beta_0=\sum_{j=1}^p\gamma_j(-t_0)^j, \qquad \beta_m=\sum_{j=m}^{p}{j \choose m}(-t_0)^{j-m}, m=1,\ldots,p,
  $$
  then the corresponding infinitesimal moments for the process $\{X(t); t\geq t_0\}$ given by $X(t)=Y(t-t_0)$ are
  $$
  A_1^X(x,t)=A_1^Y(x,t-t_0)=\eta\widetilde{P}_{\bm\gamma}(t-t_0)\,e^{-\widetilde{Q}_{\bm\gamma}(t-t_0)}=
  \alpha\,P_{\bm\beta}(t)e^{-Q_{\bm\beta}(t)}, \qquad A_2^X(x)=A_2^Y(x)
  $$
  with $\alpha=\eta\,e^{-\beta_0}$ and $P_{\bm\beta}$ the derivative of polynomial $Q_{\bm\beta}$.

  Note that $X(t)$ is also a multi-sigmoidal Gompertz diffusion process whose infinitesimal moments differ from those of $Y(t)$ in the reparametrization occurred in $\eta$ and $\bm\gamma$. The same happens for the finite dimensional distributions, transition distributions and main characteristics of the process. In particular, $E[X(t)]=E[Y(t-t_0)]$.

\section{Simulations}

In this section, some simulation examples will be carried out with the aim of illustrating the developments previously established, focusing on the strategies related to the estimation of the parameters of the model as well as the selection of the model that best fits the data. All the simulations were performed according to the following common pattern: 25 sample paths were simulated, each one obtained from expression \eqref{SolEDE}, which relates the Gompertz process under consideration and the Wiener process. All of them contain the same number of data (501), being $(i-1)\cdot 0.1$, $i=1,\ldots,501$ the observation time instants. For simplicity we have chosen a degenerate initial distribution ($P[X(0)=5]=1$). After obtaining each trajectory, we chose 51 values starting from the first one and using a step equal to 1. Hence, a sample of 51 data was obtained for each sample path.

With regard to the processes chosen for the simulation, two have been selected that correspond to situations in which there are two inflection points. The former presents a strictly increasing mean, while in the second an initial decrease is observed.

\vskip 0.25cm
\noindent\textbf{The case of increasing mean}
\vskip 0.25cm

\begin{figure}[h!]
\caption{\small Example 1. Simulated sample-paths. Black lines represent the sample mean.}
$$
\begin{array}{cc}
\includegraphics[height=5.5cm, width=7.5cm]{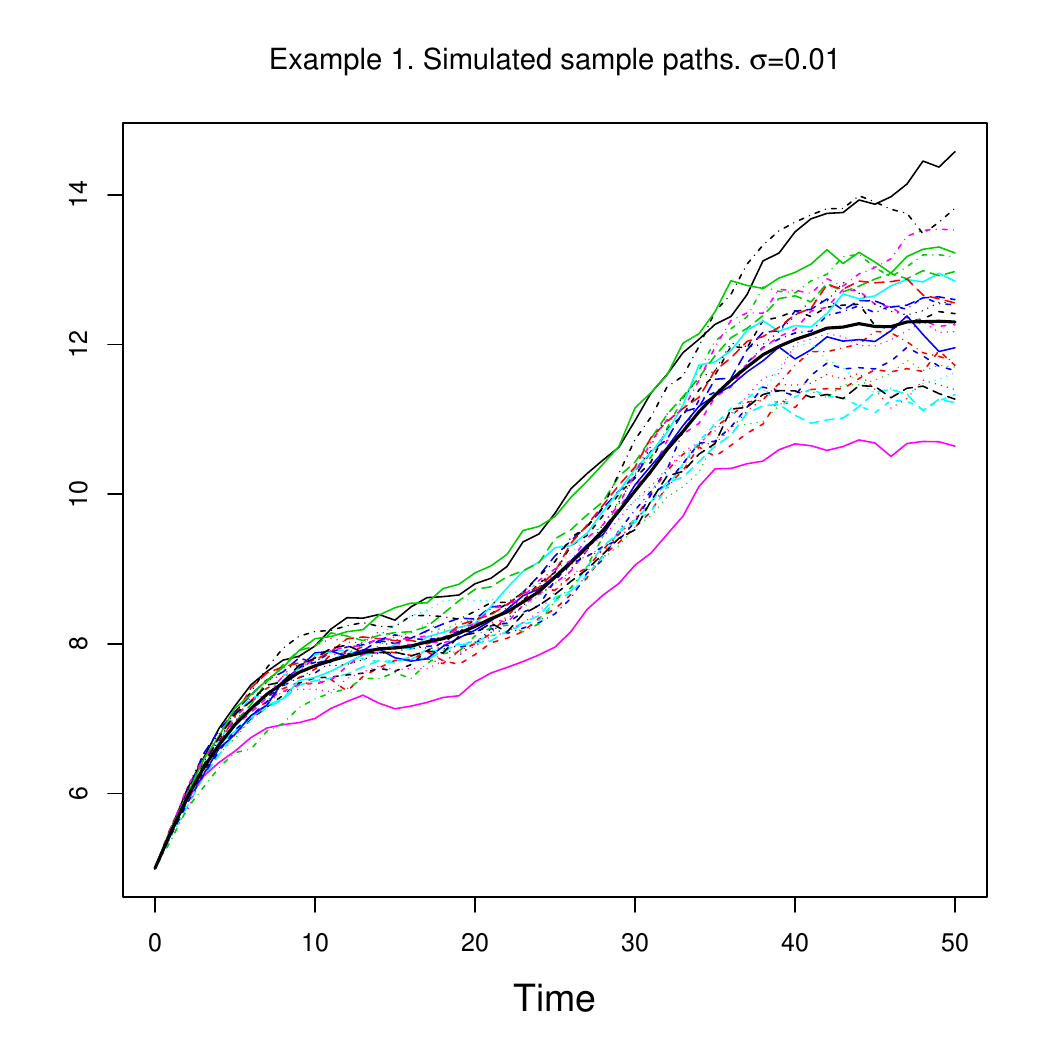} & \includegraphics[height=5.5cm,width=7.5cm]{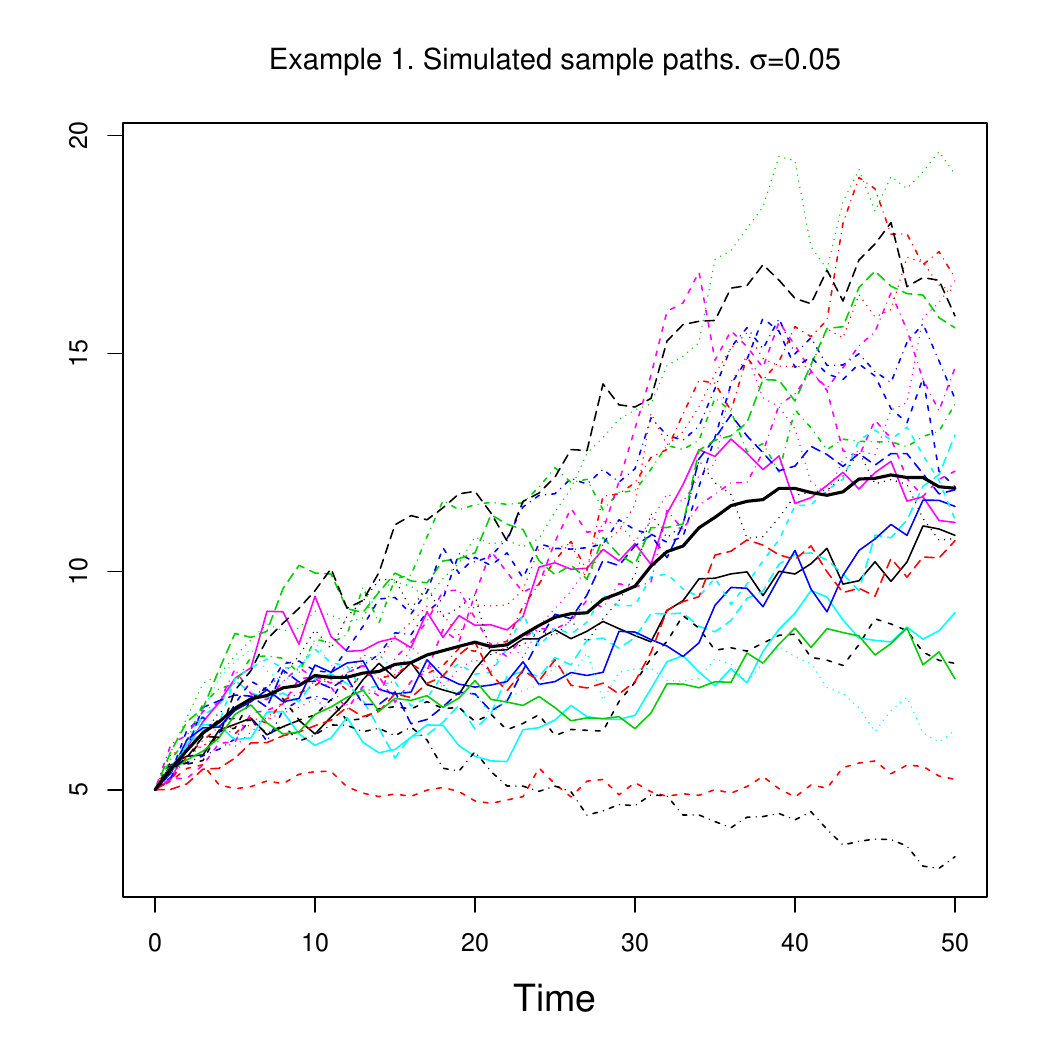}
\end{array}
$$
\label{Ejem1_1}
\end{figure}

As a first example we have selected a multi-sigmoidal Gompertz diffusion process for which the degree of the polynomial included in the infinitesimal mean is $p=3$, being $\bm\beta=(0.1225, -0.0075; 0.00017)^T$. The value of $\alpha$ is $\alpha=e^{-1}$, while two values of $\sigma$ have been considered (concretely $\sigma=0.01,0.05$) to verify the effect of increasing the infinitesimal variance in the estimation process. Figure \ref{Ejem1_1} shows the 25 simulated sample paths for each value of $\sigma$.

In order to find the process that best fits the data, for each sigma value, multi-sigmoidal Gompertz processes including polynomials from grade 2 to 5 have been considered successively. Table \ref{Ejem_1_Estimaciones} includes, for each model, the initial values of the parameters as well as their definitive estimates. The initial values have been obtained following the procedure described above. Note that the initial $\sigma$ value is common for all cases since it is calculated directly from the sample data.

\begin{table}[h!]
\centering
{\small
\begin{tabular}{c l r r r r r r r}
\hline
\multicolumn{9}{c} {$\sigma=0.01$}\\
\hline
Degree & &\multicolumn{1}{c}{$\alpha$}& \multicolumn{1}{c}{$\sigma$} & \multicolumn{1}{c}{$\beta_1$} & \multicolumn{1}{c}{$\beta_2$} & \multicolumn{1}{c}{$\beta_3$} & \multicolumn{1}{c}{$\beta_4$} & \multicolumn{1}{c}{$\beta_5$}\\
\hline
\multirow{2}{*}{2} & Initial  & 0.43264 & 0.00951 &-0.08412 & 0.00405 & & & \\
                   & Final & 0.44732 & 0.01632 & 0.19938 & 0.00726 & & & \\
\multirow{2}{*}{3} & Initial  & 0.82721 & 0.00951 & 0.09462 &-0.00576 & 0.00014 & &\\
                   & Final & 0.90536 & 0.00975 & 0.12100 &-0.00738 & 0.00016 & &\\
\multirow{2}{*}{4} & Initial  & 1.01846 & 0.00951 & 0.19476 &-0.01580 & 0.00048 &-0.0000037 &\\
                   & Final & 0.52895 & 0.01221 & 0.33611 &-0.02349 & 0.00043 & 0.0000019 &\\
\multirow{2}{*}{5} & Initial  & 0.80607 & 0.00951 & 0.01539 & 0.01266 &-0.00119 & 0.000037  & -0.00000037\\
                   & Final & 0.80607 & 0.00951 & 0.01539 & 0.01266 &-0.00119 & 0.000037  & -0.00000037\\
\hline
\multicolumn{9}{c} {$\sigma=0.05$}\\
\hline
\multirow{2}{*}{2} & Initial  & 0.40172 & 0.05147 &-0.09713 & 0.00468 & & & \\
                   & Final &-0.04547 & 0.05450 &-0.98848 & 0.28151 & & & \\
\multirow{2}{*}{3} & Initial  & 0.89617 & 0.05147 & 0.14054 &-0.00929 & 0.00021 & &\\
                   & Final & 0.87641 & 0.04832 & 0.11912 &-0.00698 & 0.00015 & &\\
\multirow{2}{*}{4} & Initial  & 0.93752 & 0.05147 & 0.16394 &-0.01181 & 0.00030 &-0.0000010 &\\
                   & Final & 0.95883 & 0.05139 & 0.15846 &-0.01143 & 0.00030 &-0.0000016 &\\
\multirow{2}{*}{5} & Initial  & 0.76767 & 0.05147 &-0.00203 & 0.01645 &-0.00147 & 0.0000458 &-0.00000043\\
                   & Final & 0.76767 & 0.05147 &-0.00203 & 0.01645 &-0.00147 & 0.0000458 & -0.00000043\\
\hline
\end{tabular}
}
\caption{Example 1. Estimates of the parameters of the models for each value of $p$.}
\label{Ejem_1_Estimaciones}
\end{table}

\begin{figure}[h!]
$$
\begin{array}{cc}
\includegraphics[height=5.5cm, width=7.5cm]{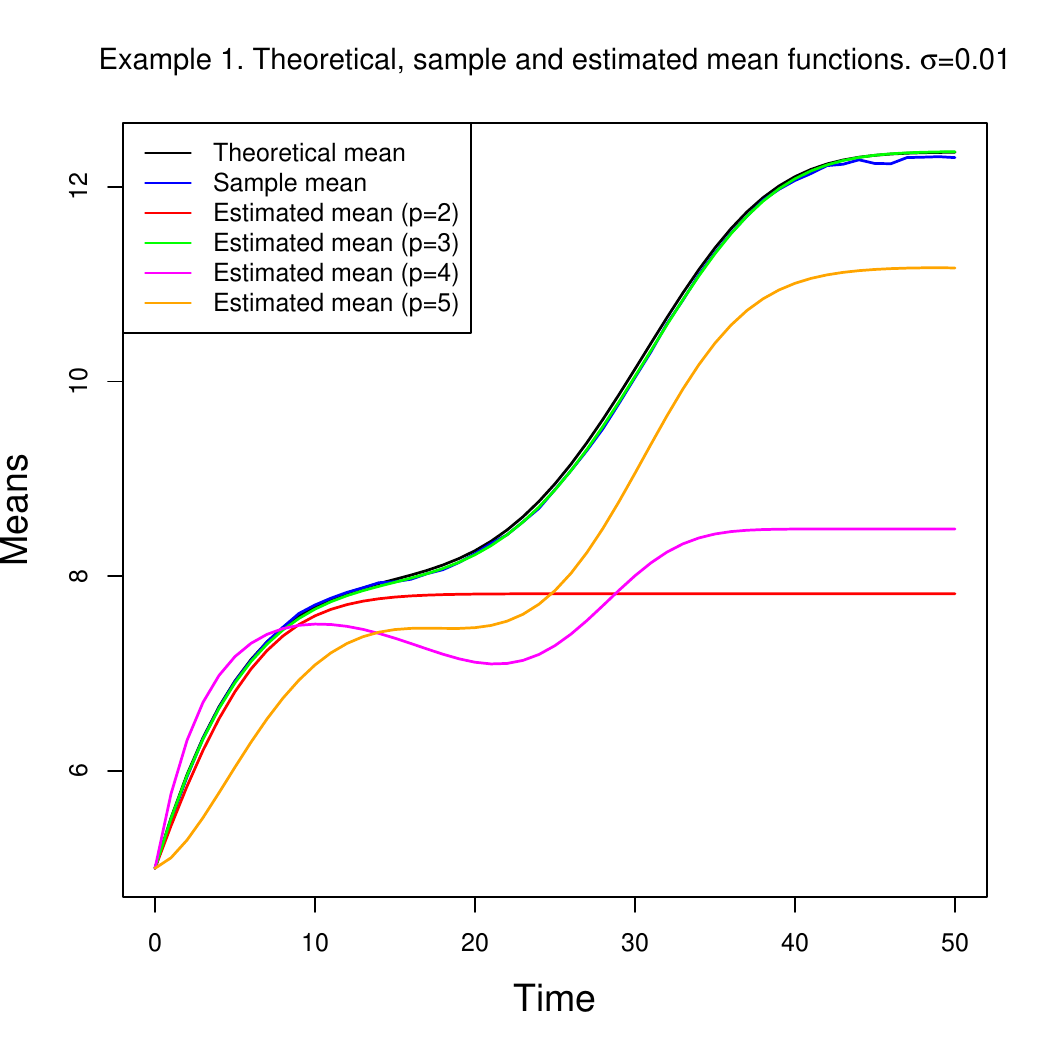} &  \includegraphics[height=5.5cm, width=7.5cm]{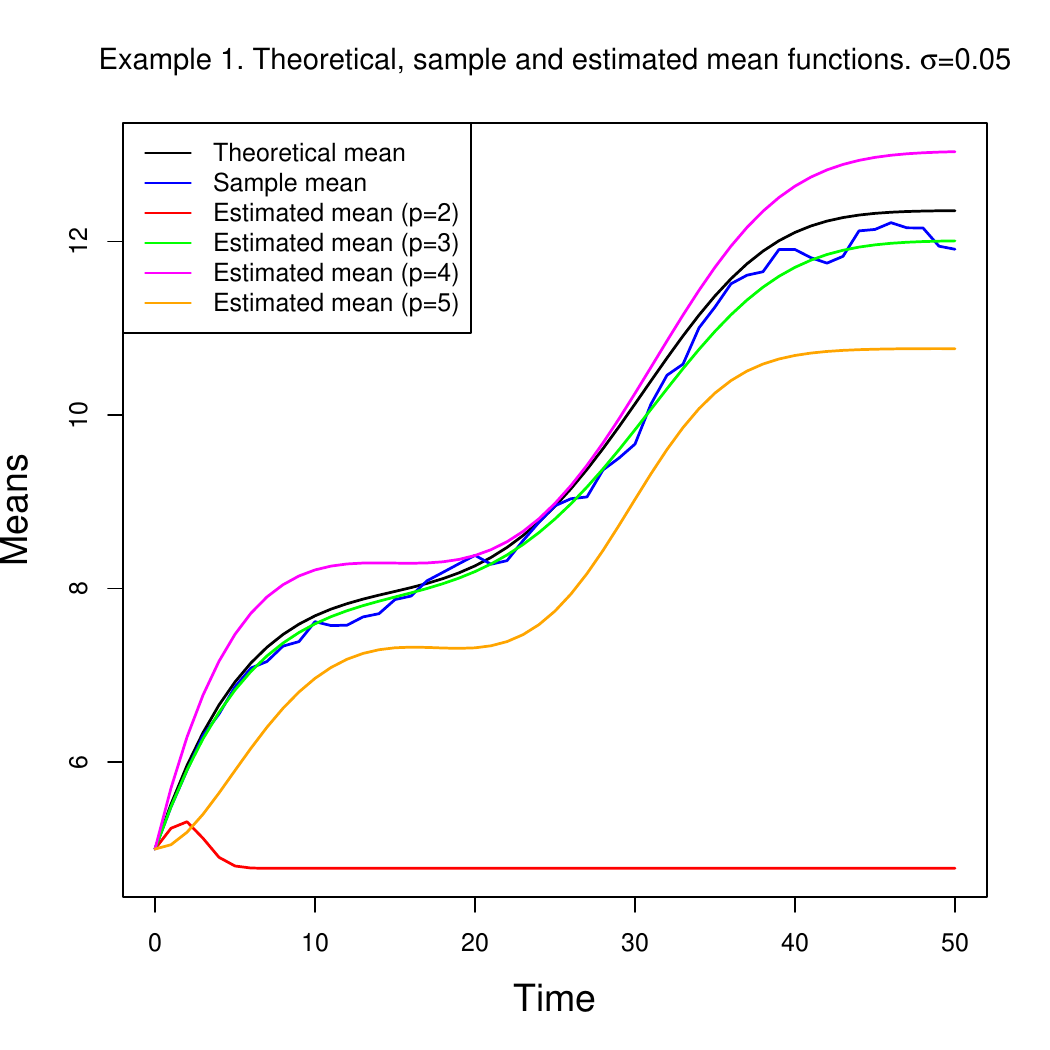}
\end{array}
$$
\caption{\small Example 1. Theoretical, sample and estimated mean functions.}
\label{Ejem1_2}
\end{figure}

Figure \ref{Ejem1_2} displays the theoretical and the sample mean functions, together with those corresponding to each estimated model. This figure suggests considering the model with $p=3$ as the optimum, although this must be endorsed by numerical measures of goodness of fit.

Table \ref{Ejem_1_Medidas} summarizes the measures that have been used (RAE, AIC and BIC). For all of them, it can be observed how from $p=3$ the goodness of fit can not be improved. Therefore, the third-degree model has been chosen as the optimal one.

\begin{table}[h!]
\centering
{\small
\begin{tabular}{c c r r r r}
\cline{3-6}
\multicolumn{2}{c}{} & \multicolumn{4}{c} {Degree}\\
\cline{3-6}
& &\multicolumn{1}{c}{2}& \multicolumn{1}{c}{3} & \multicolumn{1}{c}{4} & \multicolumn{1}{c}{5}\\
\hline
\multirow{3}{*}{$\sigma=0.01$} & RAE &  0.168  &    0.002 &   0.171 &  0.092\\
                               & AIC & -9034.903 & -10314.497 & -8952.280 & -7870.896\\
                               & BIC & -9014.380 & -10288.842 & -8921.495 & -7834.980\\
\hline
\multirow{3}{*}{$\sigma=0.05$} & RAE &  0.459  &  0.011   &  0.056  &  0.099\\
                               & AIC & -5967.582 & -6315.135  & -6275.419 & -6188.681\\
                               & BIC & -5947.058 & -6289.480  & -6244.633 & -6152.765\\
\hline
\end{tabular}
}
\caption{Example 1. Measures for choosing the estimated model.}
\label{Ejem_1_Medidas}
\end{table}

\begin{figure}[h!]
$$
\begin{array}{cc}
\includegraphics[height=5.5cm, width=7.5cm]{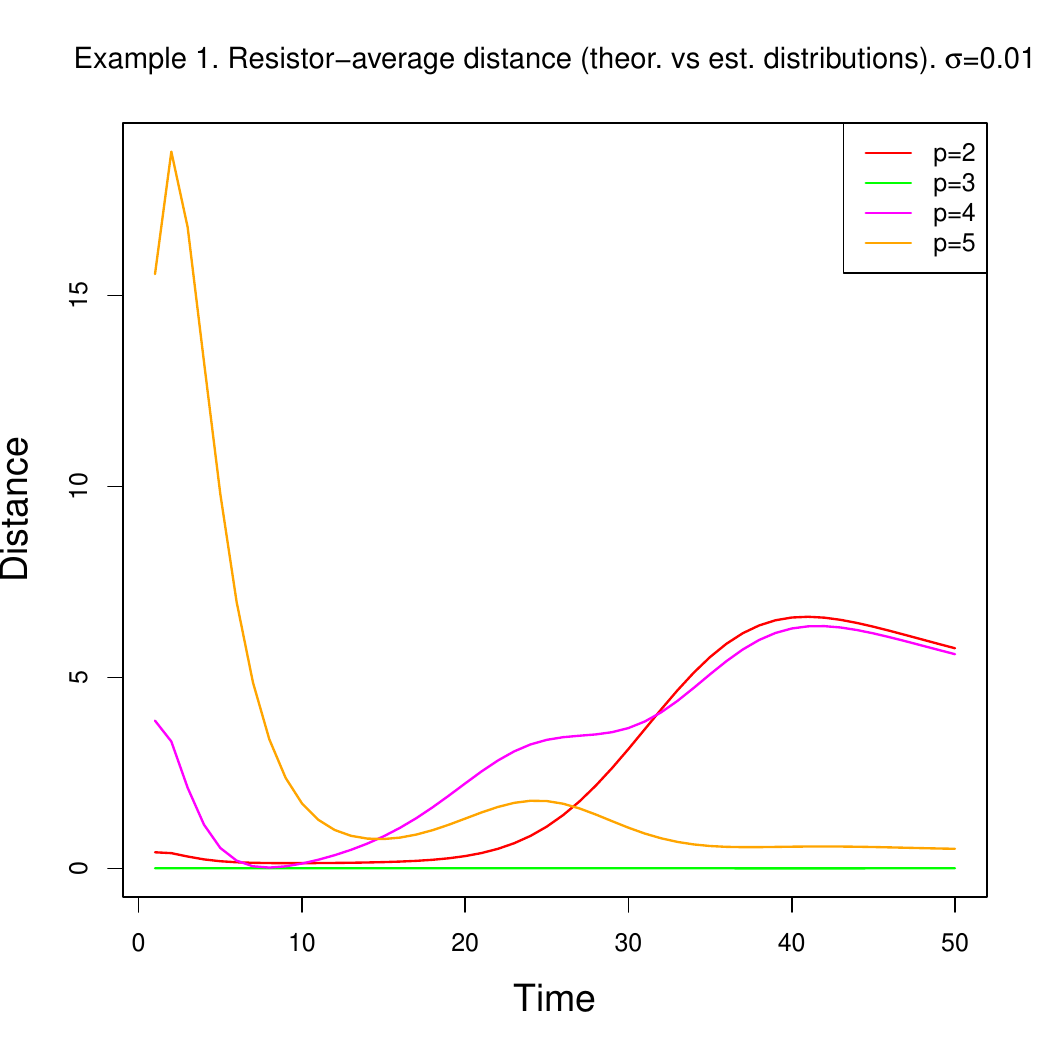} &  \includegraphics[height=5.5cm, width=7.5cm]{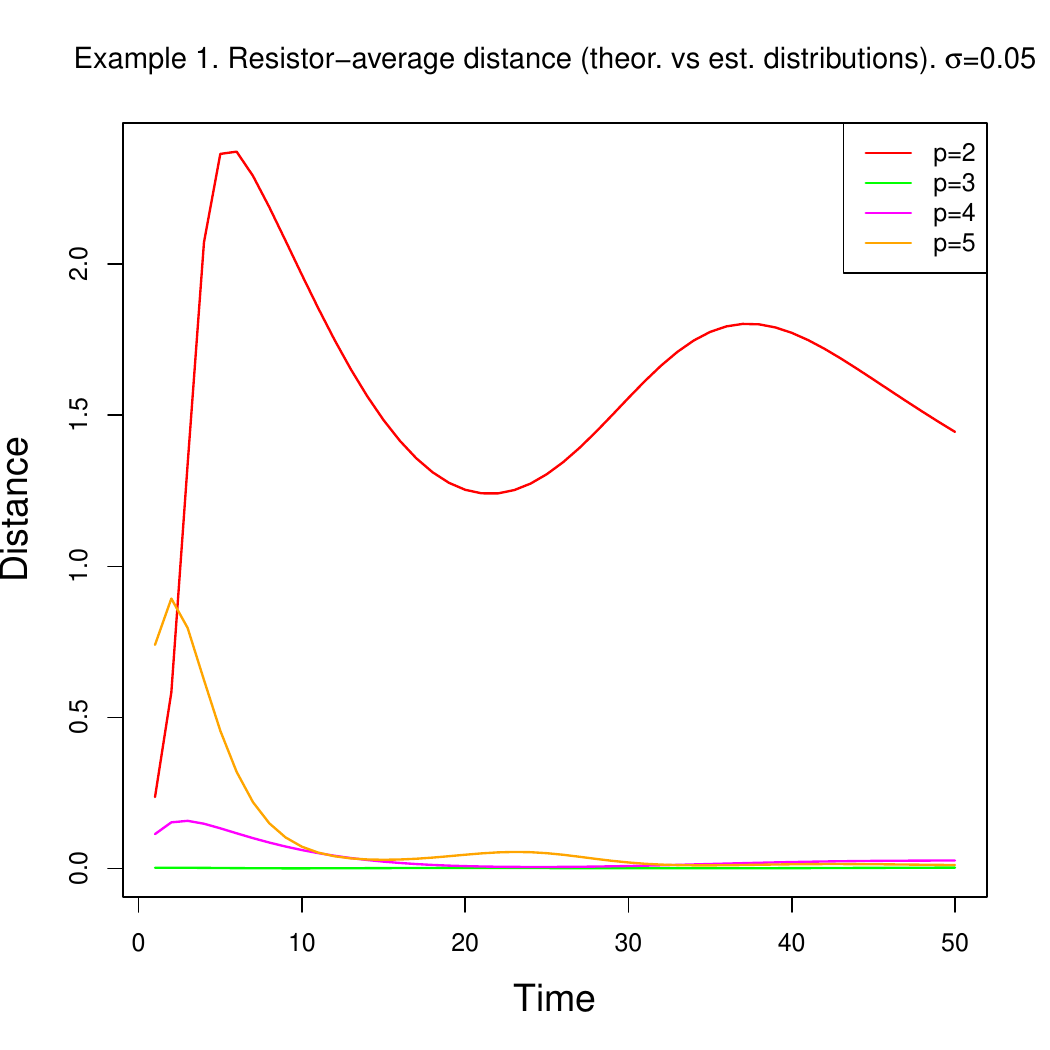}\\
\includegraphics[height=5.5cm, width=7.5cm]{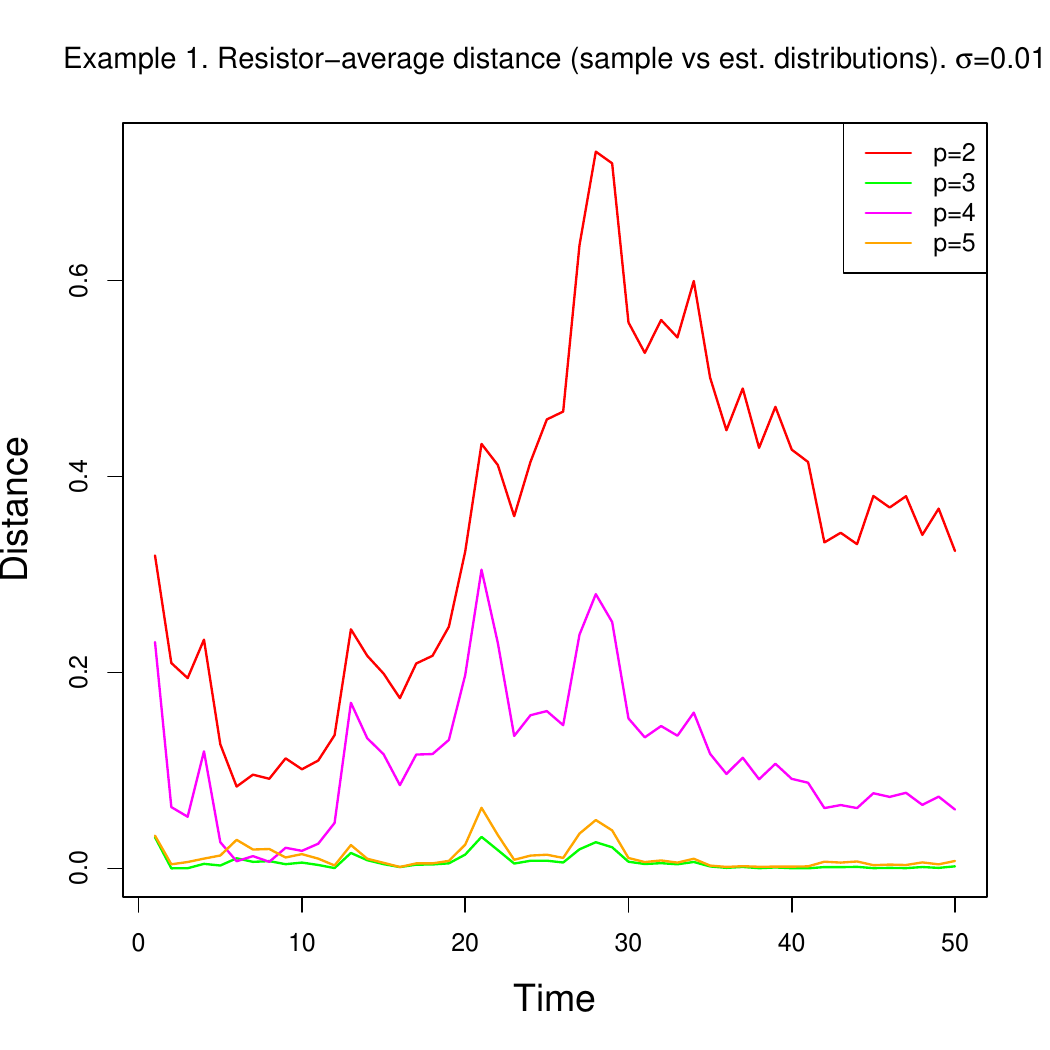} &  \includegraphics[height=5.5cm, width=7.5cm]{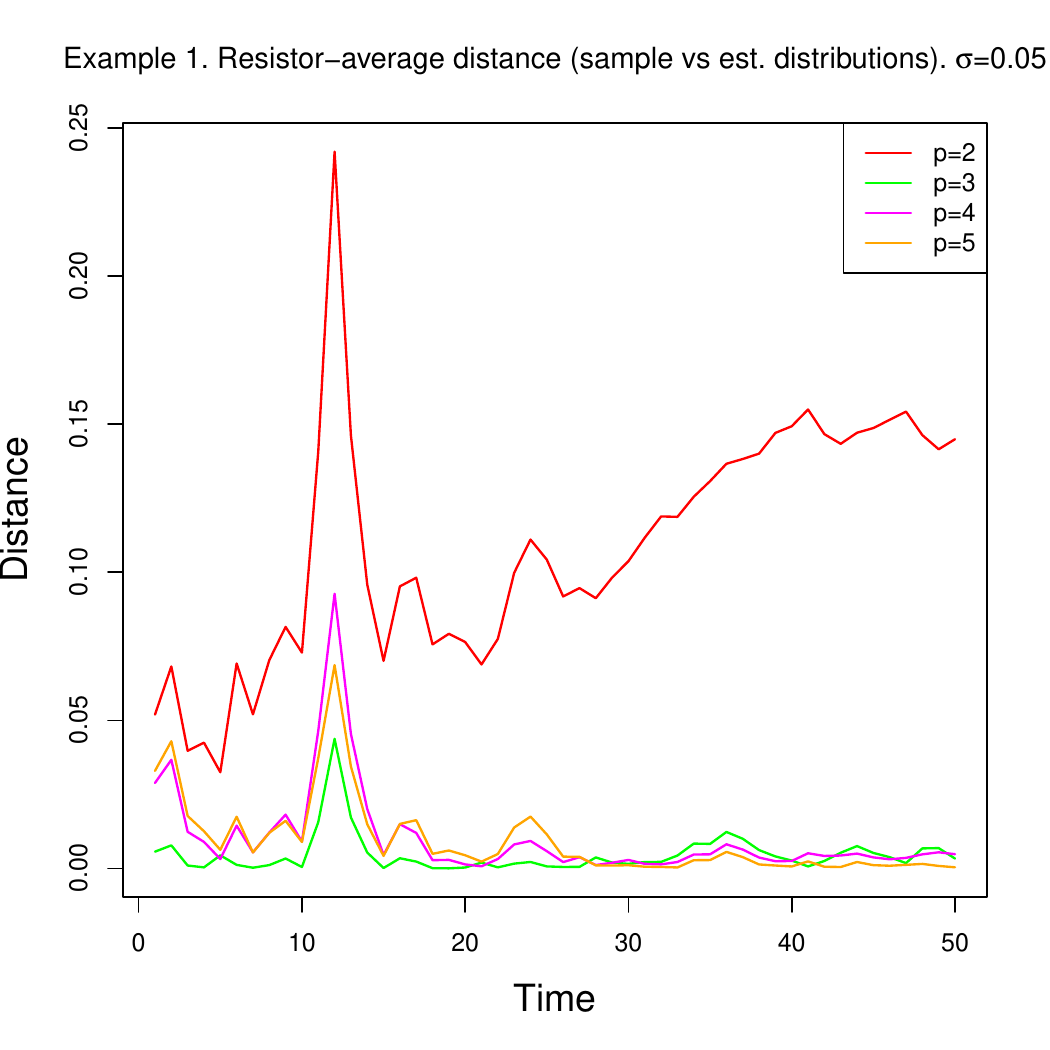}
\end{array}
$$
\caption{\small Example 1. Resistor-average distances between the theoretical and estimated models (up) and between the sample and estimated models (down).}
\label{Ejem1_3}
\end{figure}

The use of the resistor-average distance also leads us to this conclusion. For each value of $p$, and for each value of $t$, the distance between the estimated one-dimensional distribution and the corresponding theoretical and sample distributions has been calculated. This provides, for each degree, two functions whose graphs are shown in Figure \ref{Ejem1_3}, showing how odd-grade models seem to be preferable. With the idea of obtaining a globalizing measure that allows selecting the best model, Table \ref{Ejem_1_Distancias} shows, for each of them, the means and medians of the values of the distances. These two measures confirm that the model with $p=3$ is the one that should be selected as optimal. It should be noted that in practical applications, as the theoretical model is not available, the distance to be considered is that which takes the sample distribution as a reference. However, in this first example we have included the two possibilities.

\begin{table}[h!]
\centering
{\small
\begin{tabular}{c c c c c c c c c c}
\cline{3-10}
\multicolumn{2}{c}{} & \multicolumn{8}{c} {Degree}\\
\cline{3-10}
\multicolumn{2}{c}{} & \multicolumn{2}{c}{2} & \multicolumn{2}{c}{3} & \multicolumn{2}{c}{4} &
\multicolumn{2}{c}{5}\\
\cline{2-10}
& Distr. & Mean & Median & Mean & Median & Mean & Median & Mean & Median\\
\hline
\multirow{2}{*}{$\sigma=0.01$} & Theoretical & 2.71439  & 1.24435 &  0.00210 & 0.00208&  3.33965 & 3.45329 & 2.60928 & 0.95538\\
                               & Sample      & 0.34820 & 0.35109 & 0.00641 & 0.00421 & 0.11279 & 0.10989 & 0.01249 & 0.00759\\
\hline
\multirow{2}{*}{$\sigma=0.05$} & Theoretical & 1.58825 & 1.59751 & 0.00148 & 0.00139 & 0.03698 & 0.02188 & 0.10850 & 0.03090\\
                               & Sample      & 0.10873 & 0.10398 & 0.00464 & 0.00257 & 0.01023 & 0.00473 & 0.00938 & 0.00408\\
\hline
\end{tabular}
}
\caption{Example 1. Means and medians of the resistor-average distances. }
\label{Ejem_1_Distancias}
\end{table}

\begin{figure}[h!]
$$
\begin{array}{cc}
\includegraphics[height=5.5cm, width=7.5cm]{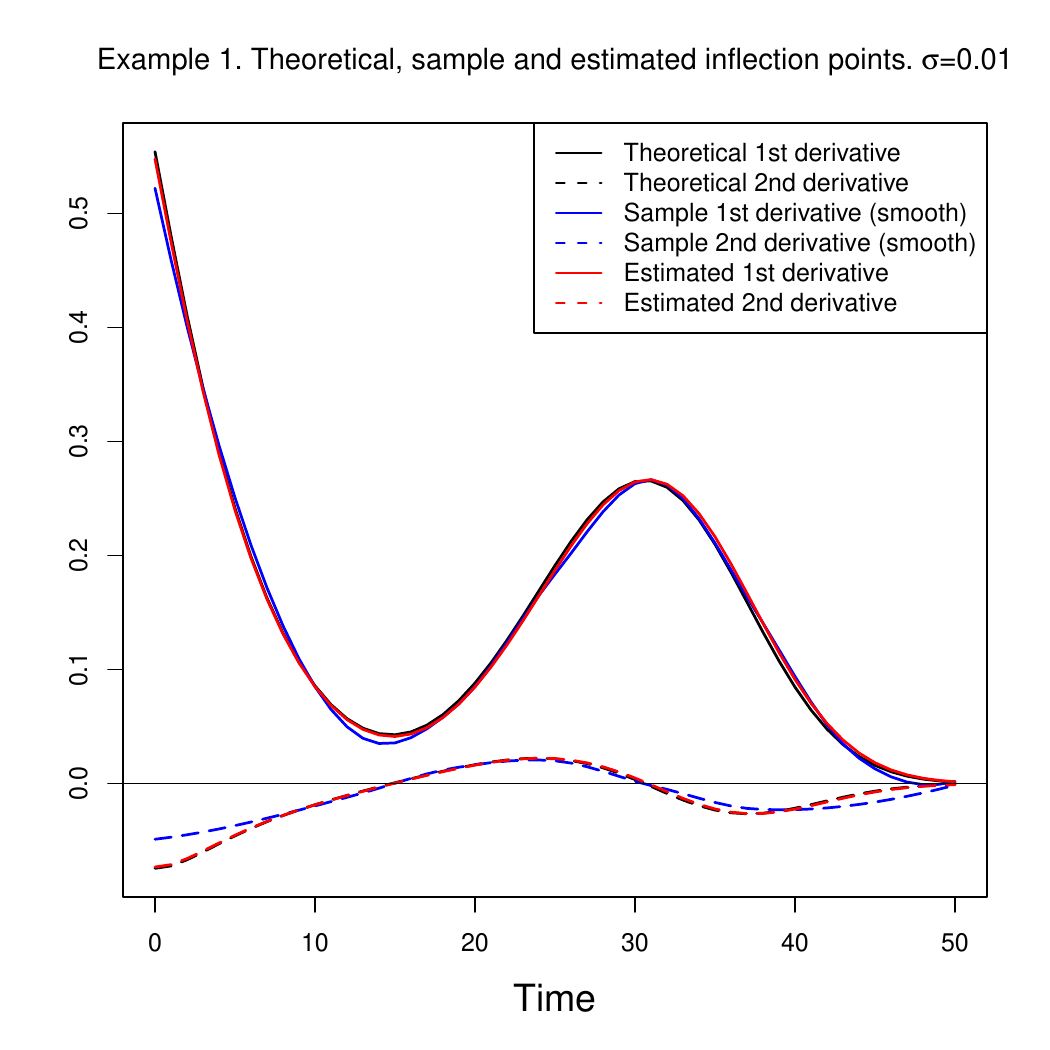} &  \includegraphics[height=5.5cm, width=7.5cm]{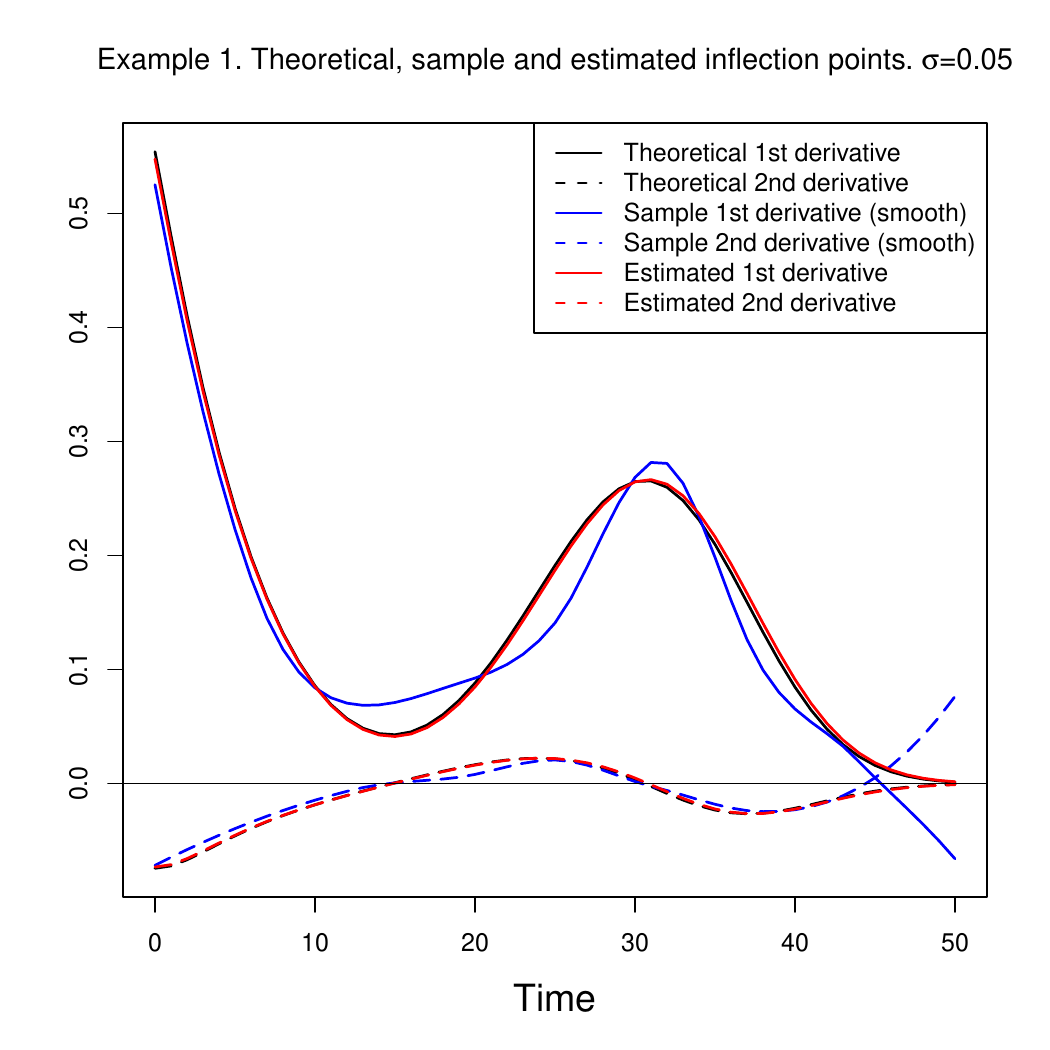}
\end{array}
$$
\caption{\small Example 1. Theoretical, sample, and estimated inflection time instants.}
\label{Ejem1_4}
\end{figure}

Finally, Figure \ref{Ejem1_4} shows, for each $\sigma$ value, the first and second derivatives of the theoretical, sample and estimated mean functions for the selected model. In order to obtain the derivatives of the sample mean function, a smoothing of the function has previously been carried out using polynomial local regression. A good fit between these functions can be observed, which is corroborated by Table \ref{Ejem_1_Inflexiones}, which contains the theoretical, sample, and estimated values of the inflection time instants.

\begin{table}[h!]
\centering
{\small
\begin{tabular}{c c c c }
\cline{3-4}
\multicolumn{2}{c}{} & \multicolumn{2}{c} {Inflection time instants}\\
\cline{3-4}
& & $t_1\, (14.787)$ & $t_2\, (30.589)$\\
\hline
\multirow{2}{*}{$\sigma=0.01$} & Sample    & 14.922 & 30.522\\
                               & Estimated & 14.869 & 30.834\\
\hline
\multirow{2}{*}{$\sigma=0.05$} & Sample    & 14.595 & 30.394\\
                               & Estimated & 14.789 & 30.553\\
\hline
\end{tabular}
}
\caption{Example 1. Sample and estimated inflection time instants for the model chosen. Theoretical values are in parentheses.}
\label{Ejem_1_Inflexiones}
\end{table}

\vskip 0.25cm
\noindent\textbf{The case of mean decreasing, then increasing}
\vskip 0.25cm

\begin{figure}[h!]
\centering
\caption{\small Example 2. Simulated sample-paths. The black line represents the sample mean.}
\includegraphics[height=6cm, width=8cm]{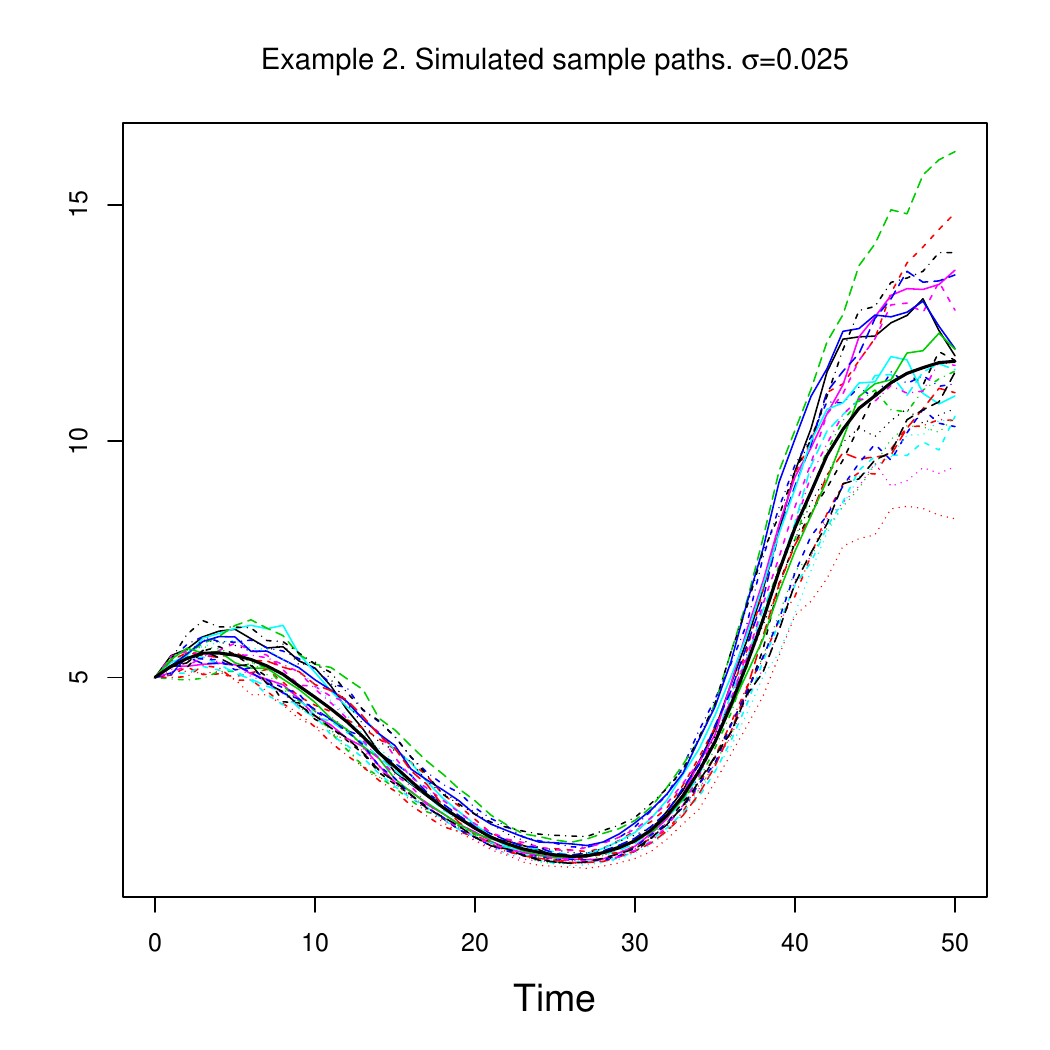}
\label{Ejem2_1}
\end{figure}

This example illustrates a case in which the data present an initial decrease and then grow up to the value of the upper bound. As in the previous case, we have selected a model with $p=3$, being $\bm\beta=(.0626,-.009,0.0002)^T$, $\alpha=e^{-1}$ and $\sigma=0.025$. Figure \ref{Ejem2_1} shows the simulated sample paths.

\begin{figure}[h!]
\centering
\includegraphics[height=6cm, width=8cm]{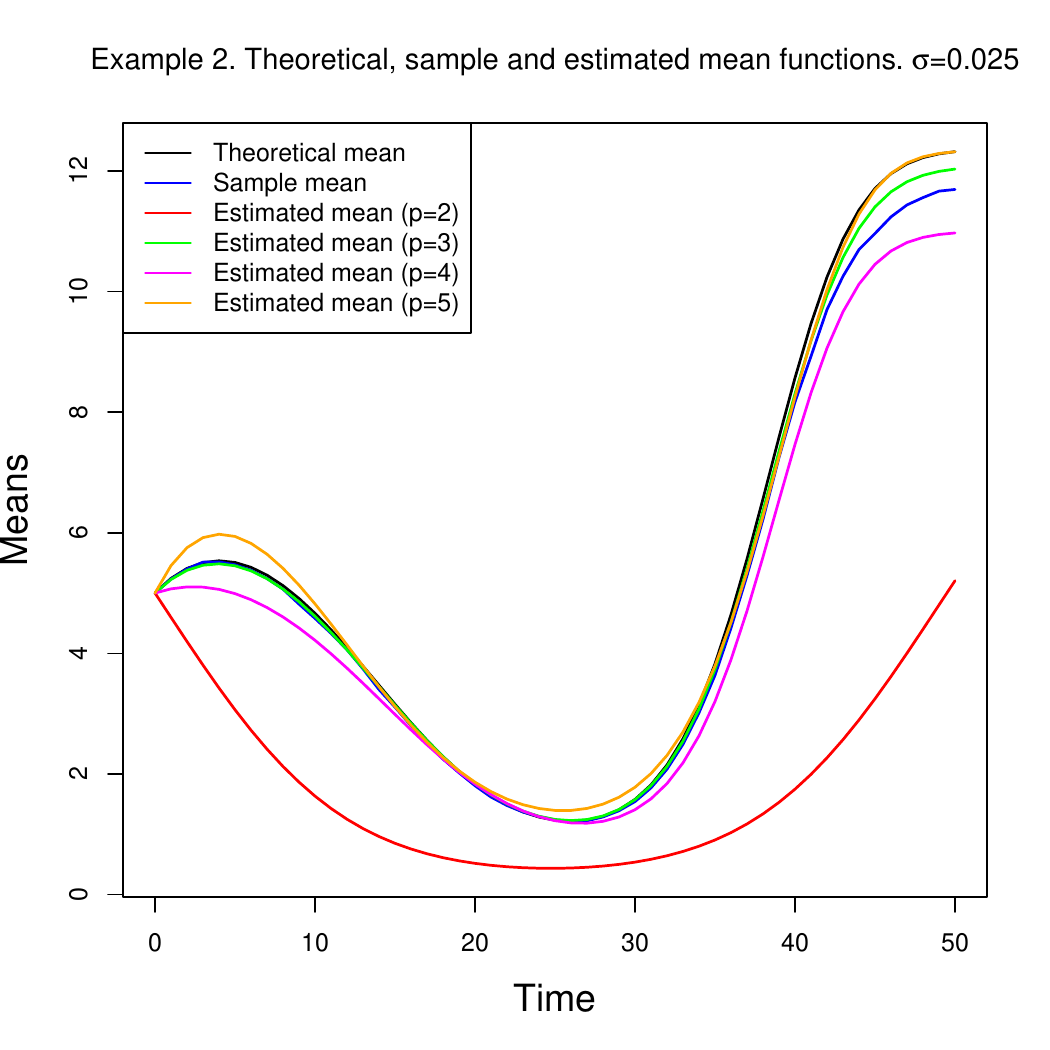}
\caption{\small Example 2. Theoretical, sample, and estimated mean functions.}
\label{Ejem2_2}
\end{figure}

Following the same methodology developed in the previous example, once again the procedure stops when considering the polynomial of degree 5. Figure \ref{Ejem2_2} shows the estimated means together with the theoretical and the sample ones.

Table \ref{Ejem_2_Medidas} contains the values of goodness-of-fit measurements, from which it follows that the model containing the third-degree polynomial must be chosen.

\begin{table}[h!]
\centering
{\small
\begin{tabular}{c r r r r}
\cline{2-5}
\multicolumn{1}{c}{} & \multicolumn{4}{c} {Degree}\\
\cline{2-5}
&\multicolumn{1}{c}{2}& \multicolumn{1}{c}{3} & \multicolumn{1}{c}{4} & \multicolumn{1}{c}{5}\\
\hline
 RAE &  0.168  &    0.002 &   0.171 &  0.092\\
 AIC & -9034.903 & -10314.497 & -8952.280 & -7870.896\\
 BIC & -9014.380 & -10288.842 & -8921.495 & -7834.980\\
\hline
\end{tabular}
}
\caption{Example 2. Measures for choosing the estimated model.}
\label{Ejem_2_Medidas}
\end{table}

\begin{figure}[h!]
$$
\begin{array}{cc}
\includegraphics[height=5.5cm, width=7.5cm]{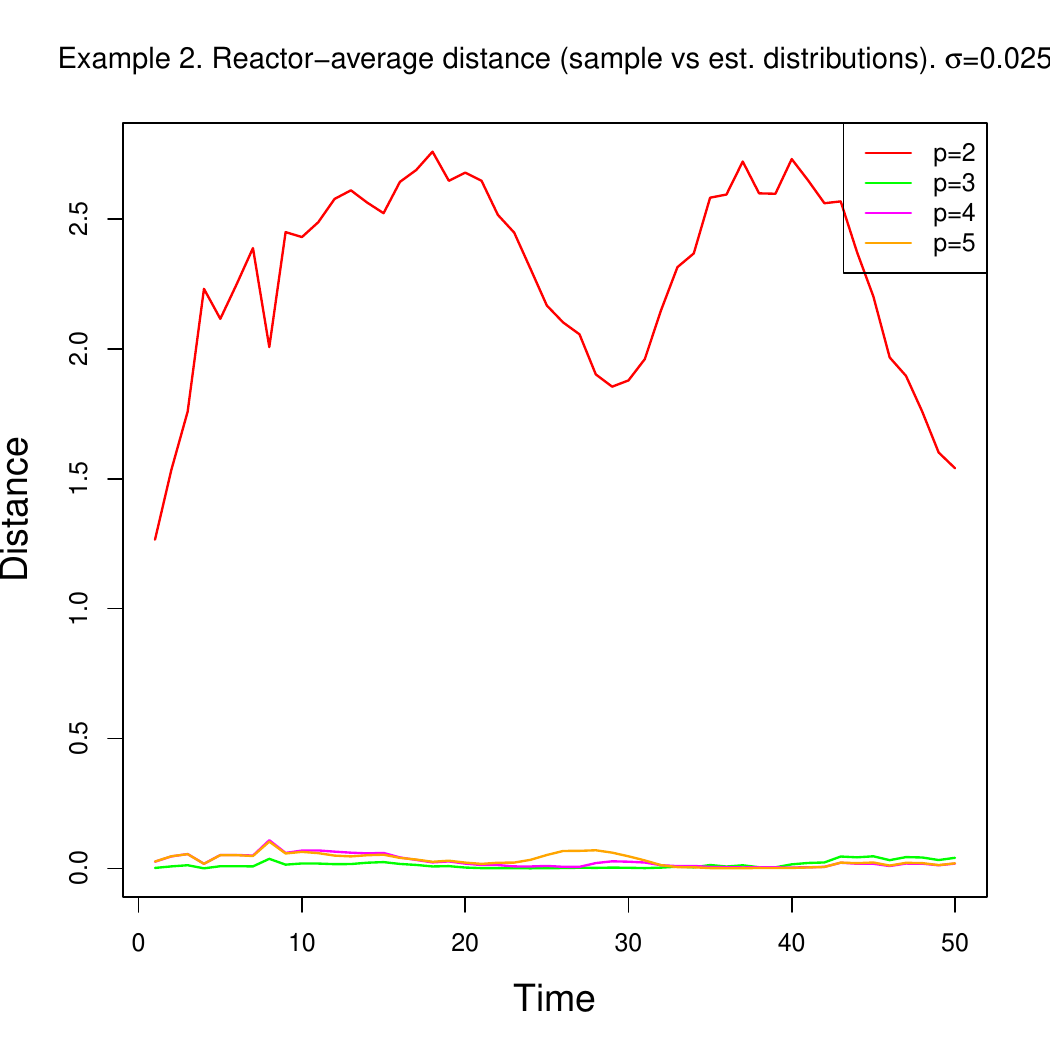} &  \includegraphics[height=5.5cm, width=7.5cm]{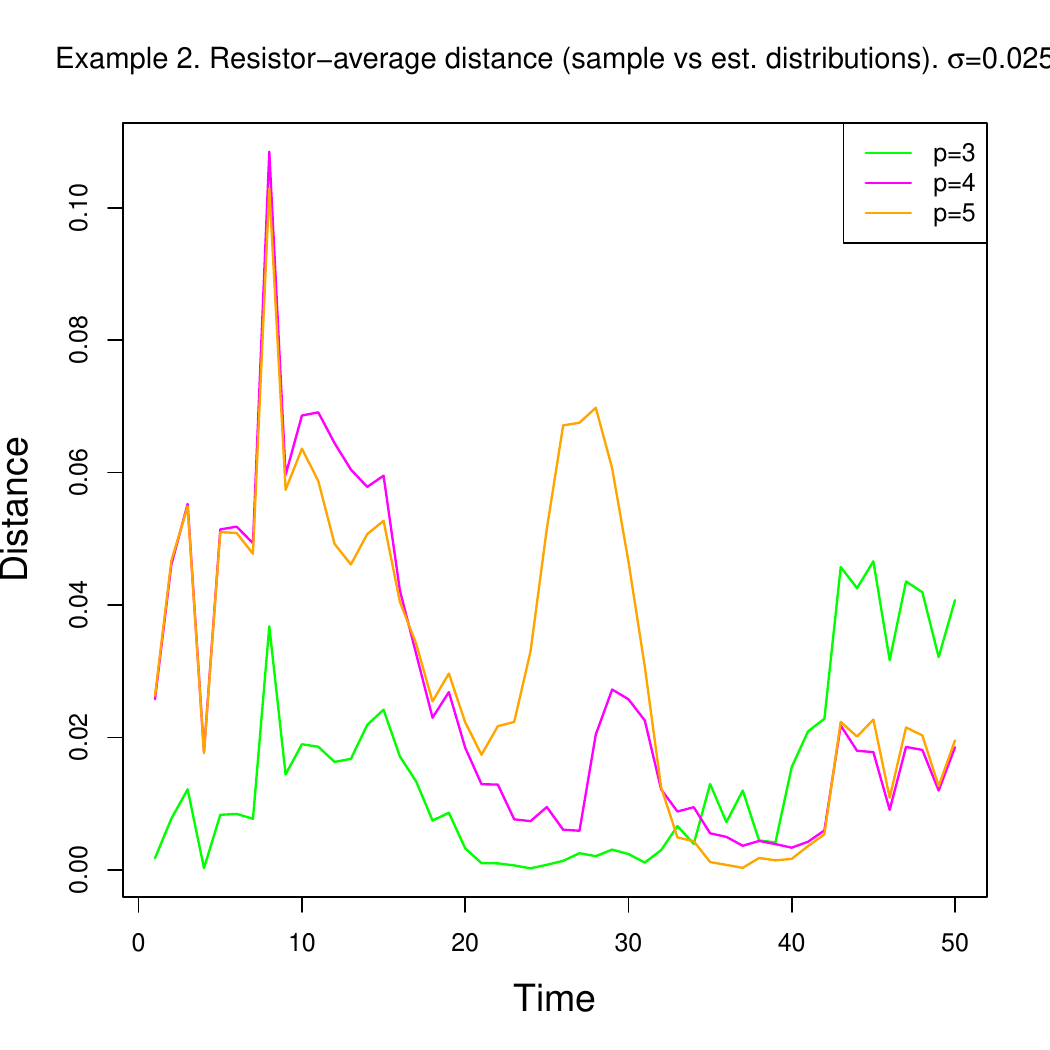}
\end{array}
$$
\caption{\small Example 2. Resistor-average distances between the sample and estimated models. The figure on the right enlarges the previous one for a better view.}
\label{Ejem2_3}
\end{figure}

Regarding resistor-average distances, in this example we have only considered those calculated between the estimated models and the sample distribution. Figure \ref{Ejem2_3} and Table \ref{Ejem_2_Distancias} contain the results obtained, confirming the previous choice of the model.

\begin{table}[h!]
\centering
{\small
\begin{tabular}{c c c c c c c c}
\hline
\multicolumn{8}{c} {Degree}\\
\hline
\multicolumn{2}{c}{2} & \multicolumn{2}{c}{3} & \multicolumn{2}{c}{4} &
\multicolumn{2}{c}{5}\\
\hline
Mean & Median & Mean & Median & Mean & Median & Mean & Median\\
\hline
3.60322 & 3.39936 & 0.00288 & 0.00270 & 0.23962 & 0.13691 & 0.15147 & 0.01742\\
\hline
\end{tabular}
}
\caption{Example 2. Means and medians of the resistor-average distances between estimated and sample unidimensional distributions.}
\label{Ejem_2_Distancias}
\end{table}

\begin{figure}[h!]
$$
\begin{array}{cc}
\includegraphics[height=6cm, width=8cm]{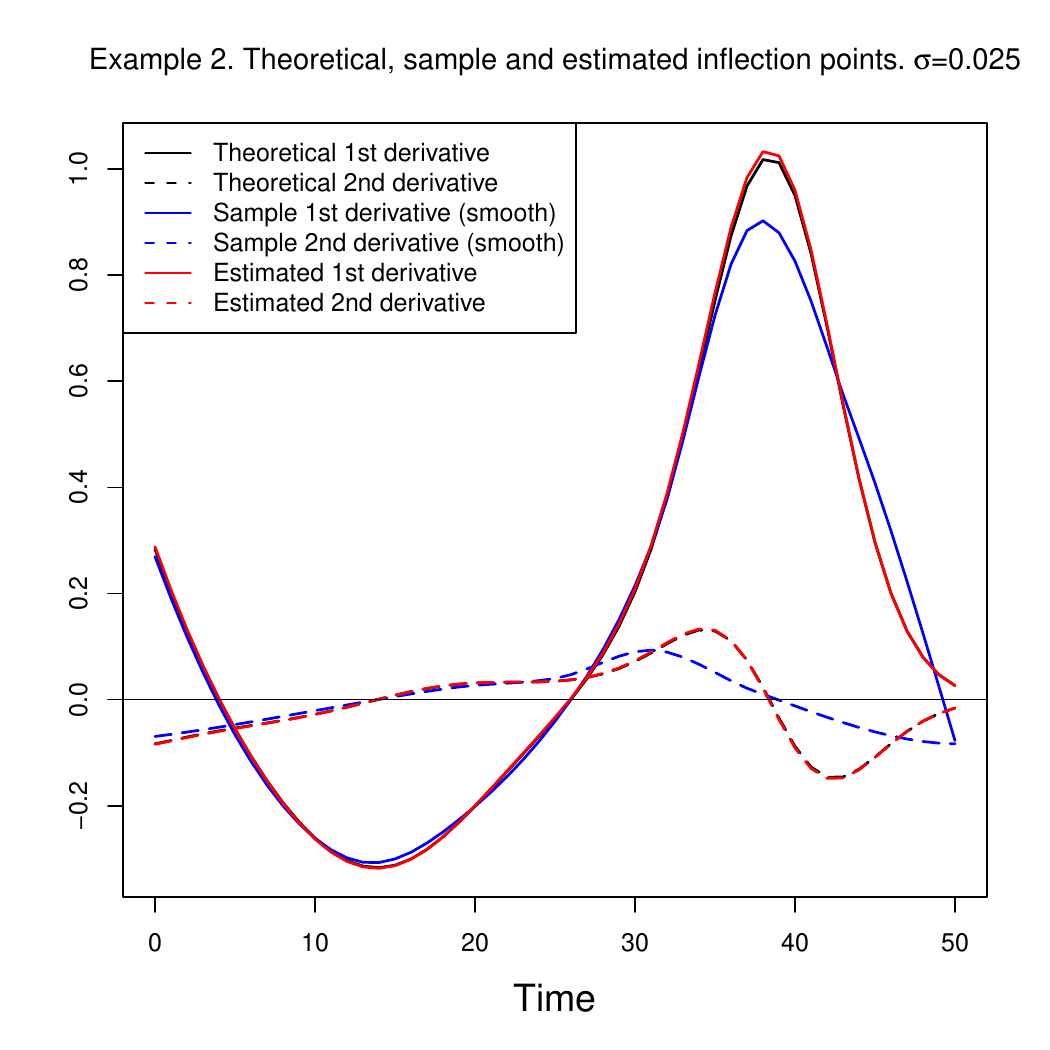}
\end{array}
$$
\caption{\small Example 2. Theoretical, sample and estimated inflection time instants.}
\label{Ejem2_4}
\end{figure}

Table \ref{Ejem_2_Inflexiones} shows the theoretical values of the inflection points (in parentheses) together with the sample and the estimated ones. In each case, the values have been obtained by numerically solving equation \eqref{DerivadaSegunda} for the theoretical, estimated, and sample mean functions. In the latter case, a natural cubic spline has been previously adjusted to the sample mean. Figure \ref{Ejem2_4} shows the situation graphically. It is readily apparent how the estimation is optimal and provides estimated values very close to the theoretical and sample values.

\begin{table}[h!]
\centering
{\small
\begin{tabular}{c c c }
\cline{2-3}
\multicolumn{1}{c}{} & \multicolumn{2}{c} {Inflection time instants}\\
\cline{2-3}
& $t_1\, (13.888)$ & $t_2\, (38.403)$\\
\hline
Sample    & 13.884 & 38.926\\
Estimated & 13.903 & 38.474\\
\hline
\end{tabular}
}
\caption{Example 2. Sample and estimated inflection time instants for the model chosen. Theoretical are in parentheses.}
\label{Ejem_2_Inflexiones}
\end{table}

The results of the simulations carried out demonstrate the suitability of the procedures developed to adjust data that follow a Gompertz multi-sigmoidal pattern. The method of obtaining initial solutions for the resolution of the system of likelihood equations provides optimal values for this purpose from the information provided by the data. The procedures introduced for calculating the degree of the polynomial operate according to two types of criteria: one based on the adjustment of the data to the model and another based on the existing proximity between the sample and estimated distributions of the process. It can be seen, in the two simulations carried out, how both types of criteria coincide in the conclusions drawn from their application.

\section{Conclusions}
A wide variety of curves (logistic, Gompertz, Bertalanffy, Richards, among others) have been used in order to describe sigmoidal growth patterns in many scientific fields. All these curves have a point of inflection that is always at a fixed proportion of its asymptotic value. Nevertheless, there are real situations in which several growth phases appear, each representing a sigmoidal pattern.

Starting from a modification of the classic Gompertz curve (by including a polynomial function in it), this paper introduces a diffusion process whose mean function is a curve of such characteristics, which allows us to model situations showing this type of behavior. The process introduced, following the methodology given in Román and Torres \cite{Rom15}, is a particular case of the lognormal process with exogenous factors, to which we apply the findings of Román et al. \cite{Rom18}. In particular, the results relative to the estimation of the process by maximum likelihood are adapted, providing strategies that provide the initial solutions required for solving the system of equations by which the parameters are estimated.

One of the main problems arising from the use of this model in real situations is that of determining the degree of the polynomial. For this purpose, the use of several criteria has been suggested. These criteria are based on the goodness of fit (relative absolute error and Akaike and BIC criteria), as well as on the measurement of the distance between the estimated and sample distributions (resistor-average distance, based on the Kullback–Leibler divergence).

\acknowledgments{This work was supported in part by the Ministerio de Econom\'ia, Industria y Competitividad, Spain, under Grant MTM2017-85568-P. The authors would also like to thank the three anonymous reviewers for their suggestions that have improved the content of the paper.}

\authorcontributions{The three authors have participated equally in the development of this work, either in the theoretical developments or in the applied aspects. The paper was also written and reviewed cooperatively.}

\conflictsofinterest{The authors declare no conflict of interest.}
\vskip 1cm

\reftitle{References}


\begin{thebibliography}{99}
\bibitem{Eby10} Eby, W.M.; Tabatabai, M.A.; Bursac, Z. Hyperbolastic modeling of tumor growth with a combined treatment of iodoacetate and dimenthylsuphoxide. \textit{BMC Cancer}. \textbf{2010}, 10:509. DOI:10.1186/1471-2407-10-509
\bibitem{Men94} Menon, A.; Mehrotra, K.; Mohan, C.K.; Ranka, S. Characterization of a class of sigmoid functions with applications to reural networks. \textit{Neural Networks}. \textbf{1996}, 9(5), 819-835. DOI:10.1016/0893-6080(95)00107-7
\bibitem{Yin03} Yin, X.; Goudriaan, J.; Lantinga, E.A.; Vos, J.; Spiertz, H. J. A flexible sigmoid function of determinate growth. \textit{Ann. Bot-London}. \textbf{2003}, 91(3), 361-371. DOI:10.1093/aob/mcg029
\bibitem{Gio07} Giovanis, A.P.: Skiadas, C.H. A new modeling approach investigating the diffusion speed of mobile telecomunication services in EU-15. \textit{Comput. Econ}. \textbf{2007}, 29(2),97–106. DOI: 10.1007/s10614-006-9067-x
\bibitem{Gal11} Gallagher, B. Peak oil analyzed with a logistic function and idealized Hubbert curve. \textit{Energ. Policy}. \textbf{2011}, 39(2), 709–802. DOI: 10.1016/j.enpol.2010.10.053
\bibitem{Yok09} Yokohama, S.; Sanada, H. Logistic regression model for predicting language change. In \textit{Issues in Quantitative Linguistics.} R. Köhler (ed). RAM-Verlag: Lüdenscheid, Germany, \textbf{2009}, 176-192. ISBN: 978-3-9802659-9-7
\bibitem{Pao12} Paolino, D.D.; Cavatorta, M.P. Sigmoidal crack growth rate: statistical modelling and applicactions. \textit{Fatig. Fract. Eng. Mater. Struct}. \textbf{2012}, 36(4), 316-326. DOI: 10.1111/ffe.12001
\bibitem{Ric79} Ricciardi, L.M. On the conjecture concerning population growth in random environment. \textit{Biol. Cybern.} \textbf{1979}, 32(2), 95–99. DOI: 10.1007/BF00337440
\bibitem{Sch07} Schurz, H. Modeling, analysis and discretization of stochastic logistic equations. \textit{Int. J. Numer. Anal. Mod.} \textbf{2007}, 4(2), 178–197.
\bibitem{Rom12} Román-Román, P.; Torres-Ruiz, F. Modelling logistic growth by a new diffusion process: Application to biological systems. \textit{Biosystems}. \textbf{2012}, 110, 9-21. DOI:10.1016/j.biosystems.2012.06.004
\bibitem{Cap74} Capocelli, R.M.; Ricciardi, L.M.  Growth with regulation in random environment. \textit{Kybernetik}. \textbf{1974}, 15(3), 147-157. DOI:10.1007/BF00274586
\bibitem{Gut07} Gutiérrez, R.; Román, P.; Romero, D.; Serrano, J.J.; Torres, F. A new gompertz-type diffusion process with application to random growth. \textit{Math. Biosci.} \textbf{2007}, 208, 147-165. DOI:10.1016/j.mbs.2006.09.020
\bibitem{Qim07} Qiming, Lv.; Pitchford, J.W. Stochastic Von Bertalanffy models, with applications to fish recruitment. \textit{J. Theor. Biol}. \textbf{2007}, 244(4), 640-655. DOI: 10.1016/j.jtbi.2006.09.009
\bibitem{Rom10} Román-Román, P.; Romero, D.; Torres-Ruiz, F. A diffusion process to model generalized von Bertalanffy growth patterns: Fitting to real data. \textit{J. Theor. Biol.} \textbf{2010}, 263(1), 59-69. DOI:10.1016/j.jtbi.2009.12.009
\bibitem{Ist17} Luz-Sant'Ana, I.; Román-Román, P.; Torres-Ruiz, F. Modeling oil production and its peak by means of a stochastic diffusion process based on the Hubbert curve. \textit{Energy}. \textbf{2017}, 133, 455-470. DOI:10.1016/j.energy.2017.05.125
\bibitem{Bar18} Barrera, A.; Román-Román, P.; Torres-Ruiz, F. A hyperbolastic type-I diffusion process: Parameter estimation bymeans of the firefly algorithm. \textit{Biosystems}. \textbf{2018}, 163, 11–22. DOI:10.1016/j.biosystems.2017.11.001
\bibitem{Tan08} Tang, S.; Heron, E.  Bayesian inference for a stochastic logistic model with switching points. \textit{Ecol. Model}. \textbf{2008}, 219, 153-169. DOI: 10.1016/j.ecolmodel.2008.08.007
\bibitem{Alb11} Albano, G.; Giorno, V.; Román-Román, P.; Torres-Ruiz, F.  Inferring the effect of therapy on tumors showing stochastic Gompertzian growth. \textit{J. Theor. Biol}. \textbf{2011}, 276, 67-77. DOI: 10.1016/j.jtbi.2011.01.040
\bibitem{Alb15} Albano, G.; Giorno, V.; Román-Román, P., Román-Román, S.; Torres-Ruiz, F. Estimating and determining the effect of a therapy on tumor dymamics by means of a modified Gompertz diffusion process. \textit{J. Theor. Biol.} \textbf{2015}, 364, 206–219. DOI: 10.1016/j.jtbi.2014.09.014
\bibitem{Rom16} Román-Román, P.; Román-Román, S.; Serrano-Pérez, J.J.; Torres-Ruiz, F.  Modeling tumor growth in the presence of a therapy with an effect on rate growth and variability by means of a modified Gompertz diffusion process. \textit{J. Theor. Biol.} \textbf{2016}, 407, 1–17. 10.1016/j.jtbi.2016.07.023
\bibitem{Alv99} Álvarez, O.; Boché, S. Modelos matemáticos para describir crecimiento doble-sigmoideos en frutos de un nectarín tardío (cv. Sun Grand). \textit{AgroSur}. \textbf{1999}, 27, 21-28.
\bibitem{Cai08} Cairns, S.P.; Robinson, D.M.; Loiselle, D.S. Double-sigmoid model for fitting fatigue profiles in mouse fast- and slow-twitch muscle. \textit{Exp. Physiol}. \textbf{2008}, 93(7), 851-862. DOI: 10.1113/expphysiol.2007.041285
\bibitem{Amo93} Amorim, L.; Bergamin Filho, A.; Hau, B.  Analysis of progress curves of sugarcane smut on different cultivars using functions of double sigmoid pattern. \textit{Phytopathology}, \textbf{1993}, 83(9), 933-936. DOI: 10.1094/Phyto-83-933
\bibitem{Rop00} Roper, L.D. Using Sigmoid and Double-Sigmoid Functions for Earth-StatesTransitions. \textbf{2000}. Available Online: http://www.roperld.com/Science/DoubleSigmoid.pdf
\bibitem{Lip10} Lipovetsky, S. Double logistic curve in regression modeling. \textit{J. Appl. Stat.} \textbf{2010}, 37(11), 1785-1793. DOI: 10.1080/02664760903093633
\bibitem{Ric83} Ricciardi, L.M.; Sacerdote, L.; Sato, S. Diffusion approximation and first-passage-time problem for a model neuron II. Outline of a computation Method. \textit{Math. Biosci.} \textbf{1983}, 64, 29-44. DOI: 10.1016/0025-5564(83)90026-3
\bibitem{Lai65} Laird, A.K. Dynamic of tumour growth: comparison of growth rates and extrapolation of growth curve to one cell. \textit{Brit. J. Cancer}, \textbf{1965}, 19(2), 278-291.
\bibitem{Rom15} Román-Román, P.; Torres-Ruiz, F. The nonhomogeneous lognormal diffusion process as a general process to model particular types of growth patterns. In \textit{Lecture Notes of Seminario Interdisciplinare di Matematica, Vol XII}, Università degli Studi della Basilicata, Potenza, Italy, \textbf{2015}, 201-219.
\bibitem{Rom18} Román-Román, P.; Román-Román, S.; Serrano-Pérez, J.J., Torres-Ruiz, F. Some Notes about inference for the lognormal diffusion process with exogenous factors. \textit{Mathematics}, \textbf{2018}, 6, 85; DOI:10.3390/math6050085
\bibitem{Joh01}  Johnson, D.H.; Sinanovic, S. Symmetrizing the Kullback-Leibler distance. \textbf{2001}.
Available at http://www.ece.rice.edu/$\sim$dhj/resistor.pdf.
\end{thebibliography}
\end{document}